\RequirePackage{ifpdf}
\ifpdf % We are running pdfTeX in pdf mode
\documentclass[pdftex]{sigma}
\else
\documentclass{sigma}
\fi

\newcommand{\rar}{\rightarrow} \newcommand{\calg}{\mathcal}

\newcommand{\strcc}[1]{\mathcal{O}_{#1}}  
\newcommand{\dcat}[1]{\text{D}^{b}(\text{\rm Sh}_{\mathbb C} [#1])}    
   
   \newcommand{\dif}{{\mathcal {D}}\text{\rm if\/f}}
\newcommand{\cale}{\mathcal{E}} \newcommand{\difdt}{{{\mathcal {D}}\text{\rm if\/f}}^{\bullet}} 
\newcommand{\diff}{\text{\rm Dif\/f}} \newcommand{\diffdt}{{\text{\rm Dif\/f}}^{\bullet}}

\newcommand{\compl}{\mathbb C} \newcommand{\dol}[2]{{\text{\rm K}^{\bullet}}_{{#1 \text{ } #2}}} 
   
 \newcommand{\fls}{I_{\text{\rm  FLS}}} \newcommand{\fmat}{{\mathfrak{gl}}_{\infty}^{\text{\rm f\/in}}}
\newcommand{\psdif}{\Psi\text{\rm Dif\/f}} \newcommand{\dga}[1]{{\mathcal B}_{#1}}

\begin{document}

\allowdisplaybreaks

\renewcommand{\PaperNumber}{010}

\FirstPageHeading

\ShortArticleName{On Lefschetz Number Formulae for Dif\/ferential Operators}

\ArticleName{Integration of Cocycles and Lefschetz Number\\ Formulae for Dif\/ferential Operators}

\Author{Ajay C.~RAMADOSS}

\AuthorNameForHeading{A.C.~Ramadoss}

\Address{Department Mathematik, ETH Z\"{u}rich, R\"{a}mistrasse 101, 8092 Z\"{u}rich, Switzerland}
\Email{\href{mailto:ajay.ramadoss@math.ethz.ch}{ajay.ramadoss@math.ethz.ch}}

\ArticleDates{Received August 12, 2010, in f\/inal form January 07, 2011;  Published online January 18, 2011}

\Abstract{Let $\cale$ be a holomorphic vector bundle on a complex manifold $X$ such that $\dim_{\compl}X=n$. Given any
continuous, basic Hochschild $2n$-cocycle $\psi_{2n}$ of the algebra $\diff_n$ of formal holomorphic dif\/ferential operators, one obtains a
$2n$-form $f_{\cale,\psi_{2n}}(\mathcal D)$ from any holomorphic dif\/ferential operator $\mathcal D$ on~$\cale$.
We apply
our earlier results
[{\em  J.~Noncommut. Geom.} {\bf 2} (2008),  405--448; {\em  J.~Noncommut. Geom.} {\bf 3} (2009), 27--45]
to show that $\int_X f_{\cale,\psi_{2n}}(\mathcal D)$ gives the Lefschetz number of
$\mathcal D$ upto a constant independent of~$X$ and $\cale$. In addition, we obtain a~``local'' result generalizing the above statement. When
$\psi_{2n}$ is the cocycle from [{\em Duke Math.~J.} {\bf 127} (2005),  487--517],
we obtain a new proof as well as a generalization of the Lefschetz number theorem of
Engeli--Felder.
We also obtain an analogous ``local'' result pertaining to B.~Shoikhet's construction of the holomorphic noncommutative residue of a~dif\/ferential operator for trivial vector bundles on complex parallelizable manifolds. This enables us to give a rigorous construction of the
holomorphic noncommutative residue of~$\mathcal D$ def\/ined by B.~Shoikhet  
when $\cale$ is an arbitrary vector bundle on an arbitrary
compact complex manifold~$X$. Our local result immediately yields a proof of a generalization of Conjecture~3.3 of~[{\em Geom. Funct. Anal.} {\bf 11} (2001), 1096--1124].}


\Classification{16E40; 32L05; 32C38; 58J42}

\section{Introduction}\label{section1}

\textbf{1.1.} Let $X$ be a connected compact complex manifold such that $\dim_{\compl}X=n$. Let $\cale$ be a~holomorphic vector bundle on
$X$. In this paper, the term ``vector bundle'' shall hence forth mean ``holomorphic vector bundle''. Let $\dif(\cale)$ denote the sheaf of
holomorphic dif\/ferential operators on $\cale$. Let $\Omega^{0,\bullet}_X$ denote the Dolbeault resolution of the sheaf $\strcc{X}$ of
holomorphic functions on $X$. Let $\difdt(\cale)=\Omega^{0,\bullet}_X \otimes_{\strcc{X}} \dif(\cale)$ and let
$\diffdt(\cale)=\Gamma(X,\difdt(\cale))$. There is a suitable topo\-lo\-gy on $\diffdt(\cale)$ which we describe in Section~\ref{section2}. Recall that any
global holomorphic dif\/ferential operator $\mathcal D$ induces an endomorphism of the Dolbeault complex $\dol{\cale}{}:=
\Gamma(X,\Omega^{0,\bullet}_X \otimes_{\strcc{X}} \cale)$ that commutes with the $\bar{\partial}$ dif\/ferential. The {\it Lefschetz number} of
$\mathcal D$, also known as the supertrace of~$\mathcal D$,
\begin{gather*}
\text{str}(\mathcal D):=\sum_i {(-1)}^i\text{Tr}(\mathcal
D)|_{\text{H}^i(X,\cale)}
\end{gather*}
 therefore makes sense. More generally, let $\cale$ be a vector bundle with bounded geometry on an arbitrary
complex manifold $X$ of complex dimension $n$. The reader may refer to Section~2.2 for the def\/inition of ``bounded geometry''. There is a
complex $\widetilde{\text{hoch}(\dif(\cale))}$ of ``completed'' Hochschild chains that one can associate with $\difdt(\cale)$. This is a complex
of (c-soft) sheaves of $\compl$-vector spaces on $X$. Let $\Gamma_c$ denote ``sections with compact support''. Suppose that $\alpha \in
\Gamma_c(X,\widetilde{\text{hoch}(\dif(\cale))})$ (see Section~\ref{section2} for precise def\/initions). Let $\alpha_0$ denote the component of $\alpha$
with antiholomorphic degree~$0$. For any $t>0$, $\alpha_0\text{e}^{-t\Delta_{\cale}}$ is a trace class operator on the (graded) Hilbert space
$\dol{\cale}{L^2}$ of square integrable elements of $\dol{\cale}{}\!\!$ (Theorem~\ref{oldrecall2}, part~1). The {\it supertrace}
$\text{str}(\alpha_0e^{-t\Delta_{\cale}})$ of $\alpha_0\text{e}^{-t\Delta_{\cale}}$ is given by the formula
\begin{gather*}
\text{str}\big(\alpha_0\text{e}^{-t\Delta_{\cale}}\big)= \sum_i {(-1)}^i
\text{Tr}\big(\alpha_0\text{e}^{-t\Delta_{\cale}}\big)|_{\text{K}_{\cale, L^2}^i}.
\end{gather*}
We may def\/ine the Lefschetz number of $\alpha$ to be
\begin{gather*}
\lim_{t \rar \infty} \text{str}\big(\alpha_0e^{-t\Delta_{\cale}}\big).
\end{gather*}
 When $X$ is compact and $\alpha=\mathcal D$ is a global
holomorphic dif\/ferential operator on $\cale$, this coincides with the Lefschetz number of $\mathcal D$ def\/ined earlier (any vector bundle on a
compact manifold is of bounded geometry).

\textbf{1.2.} Let $\diff_n$ denote the algebra $\compl[[y_1,\dots,y_n]][\partial_1,\dots,\partial_n]$ ($\partial_i:=\frac{\partial}{\partial y_i}$)
of formal holomorphic dif\/ferential operators. Denote the complex of continuous normalized Hochschild  cochains of $\diff_n$ with
coef\/f\/icients in $\diff_n^*$ by $\text{C}^{\bullet}(\diff_n)$. Let $\psi_{2n}$ be a cocycle in $\text{C}^{2n}(\diff_n)$. For an al\-gebra~$A$, let~$\text{M}_r(A)$ denote the algebra of $r \times r$ matrices with entries in~$A$. Then $\psi_{2n}$ extends to a~cocycle~$\psi^r_{2n}$ in $\text{C}^{2n}(\text{M}_r(\diff_n))$
by the formula
\begin{gather} \label{eq1} \psi^r_{2n}(A_0 \otimes D_0,\dots,A_{2n+1} \otimes
D_{2n+1})=\text{tr}(A_0\cdots A_{2n+1})\psi_{2n}(D_0,\dots,D_{2n+1}).
\end{gather} Note that the group $G:= {GL}(n,\compl) \times
{GL}(r,\compl)$ acts on the algebra $\text{M}_{r}(\diff_n)$ with $GL(r,\compl)$ acting by conjugation and $GL(n,\compl)$ acting by linear
coordinate changes. It can be checked without dif\/f\/iculty from  \eqref{eq1}  that if $\psi_{2n}$ is $GL(n)$-{\it basic}, then $\psi^r_{2n}$ is
$G$-basic. We will def\/ine the term basic carefully in Section~3.2.2.

Given a basic Hochschild $2n$-cocycle $\psi_{2n}$, a Gel'fand--Fuks type procedure (which extends a~construction from \cite{enfe} and
\cite{pfpota}) enables one to construct a map of complexes of sheaves
\begin{gather*}
f_{\cale,\psi_{2n}}: \ \ \widetilde{\text{hoch}(\dif(\cale))} \rar
\Omega^{2n-\bullet}_X,
\end{gather*}
 where the target has the usual de~Rham dif\/ferential (see Section~3.2 for details). Our extension of the construction
from \cite{enfe,pfpota} is made possible by a result from \cite{CR} (Theorem \ref{fedosov}) showing the existence of a ``good'' choice
of Fedosov connection on the bundle $\dga{\cale}$ of f\/ibrewise formal dif\/ferential operators on $\cale$. $f_{\cale,\psi_{2n}}$ induces a map
of complexes
\begin{gather*}
f_{\cale,\psi_{2n}}: \ \ \Gamma_c(X,\widetilde{\text{hoch}(\dif(\cale))}) \rar \Omega^{2n-\bullet}_c(X).
\end{gather*}
 The subscript `$c$' above
 indicates
``compact supports''. As a map of complexes of
sheaves, $f_{\cale,\psi_{2n}}$~depends upon the choice of a Fedosov connection on $\dga{\cale}$. Nevertheless,
 as a map in the bounded derived category $\dcat{X}$ of sheaves of $\compl$-vector spaces on~$X$, $f_{\cale,\psi_{2n}}$ is independent of the
 choice of Fedosov
 connection on the bundle of f\/ibrewise formal dif\/ferential operators on~$\cale$ (see Proposition~\ref{welldefined}). Let $\tau_{2n}$ denote
 the normalized
 Hochschild cocycle from \cite{fefesh}. Let $[\cdots]$ denote ``class in cohomology''. Since $\text{HH}^{2n}(\diff_n)$ is $1$-dimensional, the
 ratio
 $\frac{[\psi_{2n}]}{[\tau_{2n}]}$ makes sense.

\textbf{1.3.} One of the main results in this paper (proven in Section~\ref{section3}) is the following.

\begin{theorem} \label{th} If $\alpha \in \Gamma_c(X,\widetilde{\text{\rm hoch}(\dif(\cale))})$ is a $0$-cycle,
\begin{gather*}
\int_X f_{\cale,\psi_{2n}}(\alpha) =
(2 \pi i)^n \frac{[\psi_{2n}]}{[\tau_{2n}]} . \lim_{t\rar \infty} \text{\rm str}\big(\alpha_0\text{e}^{-t\Delta_{\cale}}\big).
\end{gather*}
\end{theorem}

A special case of the above theorem, when $\alpha=\mathcal D$ is a global holomorphic dif\/ferential operator on $X$ and when $X$ is compact is
as follows.

\begin{theorem} \label{th1}
\begin{gather*}
\int_X f_{\cale,\psi_{2n}}(\mathcal D)=(2 \pi i)^n \frac{[\psi_{2n}]}{[\tau_{2n}]}.\text{\rm str}(\mathcal D).
\end{gather*}
\end{theorem}

Theorem \ref{th} is proven using results from \cite{Ram1} and \cite{Ram2}. If $\psi_{2n}$ is known explicitly, $f_{\cale}(\alpha)$ can in
principle be locally calculated explicitly for $\alpha \in \Gamma_c(X,\widetilde{\text{hoch}(\dif(\cale))})$.

\textbf{1.4} In particular, if we put $\psi_{2n}$ to be the normalized Hochschild cocycle $\tau_{2n}$ of \cite{fefesh}, the $2n$-form
$f_{\cale,\tau_{2n}}(\mathcal D)$ was denoted in \cite{enfe} by $\chi_0(\mathcal D)$. We therefore obtain a dif\/ferent proof of the following
Lefschetz number theorem of \cite{enfe} as well as a generalization of this result.

\medskip

\noindent
\textbf{Engeli--Felder theorem \cite{enfe}.} {\it If $X$ is compact and $\mathcal D \in \Gamma(X,\dif(\cale))$, }
\begin{gather*}
\text{\rm str}(\mathcal D) =
\frac{1}{(2 \pi i)^n} \int_X \chi_0(\mathcal D).
\end{gather*}

We sketch in a remark Section~3.5 why the map $\alpha \mapsto \int_X f_{\cale,\psi_{2n}}(\alpha)$ may be viewed as an ``integral of
$\psi^r_{2n}$ over $X$'' in some sense. Upto this point, this paper completes our work in \cite{Ram1} and \cite{Ram2} by combining the results therein with an extension of the trace density construction of~\cite{enfe} using a more delicate choice of Fedosov connection on ${\mathcal B}_{\cale}$.

\textbf{1.4.1.} Let us brief\/ly summarize our approach to Lefschetz number formulae for dif\/ferential operators in \cite{Ram1,Ram2} and the current paper. Our approach is to show that the two sides of Theorem \ref{th} coincide upto a constant fact $C(\psi_{2n})$ independent of $X$ and $\cale$ such that $C(\psi_{2n})=\frac{[\psi_{2n}]}{[\tau_{2n}]}C(\tau_{2n})$. Compared to the approach in \cite{enfe}, the analytic routine part of this step is minimal: beyond proving that $\alpha \mapsto \lim\limits_{t\rar \infty} \text{str}(\alpha_0\text{e}^{-t\Delta_{\cale}})$ vanishes on $0$-boundaries (which required only the most basic heat kernel estimates (see \cite{Ram1})), this step requires a few simple, natural arguments from \cite{Ram1, Ram2} and the current paper. $C(\tau_{2n})$ can then be determined by computing both sides of Theorem \ref{th} for the identity operator on ${({\mathbb P}^1)}^{\times n}$ (see \cite{Ram2}). In particular, the local, more general result of Theorem \ref{th} is an integral part of our approach and the Engeli--Felder formula follows as a direct corollary of this result. As far as this author can see, local results like Theorem \ref{th} do not appear in the approach in \cite{enfe}.

\textbf{1.5.} We proceed to prove Conjecture~3.3 from \cite{Sh1} in greater generality. In \cite{Sh1}, B.~Shoikhet starts with a particular
$2n+1$ cocycle $\Psi_{2n+1} \in \text{C}^{2n+1}_{\text{Lie}}(\fmat(\diff_n),\compl)$ called a lifting formula. He then attempts to construct a
map of complexes
\begin{gather*}
\Theta: \ \ \text{C}^{\rm lie}_{\bullet}\big(\fmat{(\diff(\cale))},\compl\big) \rar \Omega^{2n+1-\bullet}(X).
\end{gather*} For any global
holomorphic dif\/ferential operator $\mathcal D$ on $\cale$, he refers to the number $\int_X \Theta(E_{11}(\mathcal D))$ as its {\it holomorphic
noncommutative residue} (whenever $X$ is compact). Since $\Psi_{2n+1}$ is not $GL(n)$-basic, this cannot be done as directly as in~\cite{Sh1}
unless $X$ is complex parallelizable and $\cale$ is trivial (as per our understanding).

One can construct a complex of ``completed'' Lie chains $\widetilde{\text{Lie}}(\dif(\cale))$ of $\fmat(\dif(\cale))$ (see Section~2.4). This
complex is analogous to $\widetilde{\text{hoch}}(\dif(\cale))$. For a trivial vector bundle $\cale$ (with bounded geometry) on a
parallelizable (not necessarily compact) complex manifold $X$, we construct a map of complexes of sheaves
\begin{gather*}
\lambda(\Psi^r_{2n+1}): \ \ \widetilde{\text{Lie}}(\dif(\cale))[1] \rar \Omega^{2n-\bullet}_X
\end{gather*}
 using B.~Shoikhet's lifting formulas. Let $\alpha
\in \Gamma_c(X,\widetilde{\text{Lie}}(\dif(\cale))[1])$ be a nontrivial $0$-cycle. The component $\alpha_0$ of $\alpha$ of antiholomorphic
degree $0$ is a f\/inite matrix with entries in $\diff^0(\cale)$.

\begin{theorem}\label{th3}
\begin{gather*}
\int_X \lambda(\Psi_{2n+1})(\alpha)=C.\lim_{t \rar \infty} \text{\rm str}\big(\text{\rm Tr}(\alpha_0)e^{-t\Delta_{\cale}}\big),
\end{gather*}
 where $C$ is a
constant independent of $X$ and $\cale$. \end{theorem}

We prove Theorem~\ref{th3} in Section~4.2. The methods used for this are similar to those used in~\cite{Ram1} and~\cite{Ram2}. In Sections~4.2.5 and~4.2.6, we provide a rigorous (though not as explicit as we had originally hoped for) def\/inition of B.~Shoikhet's holomorphic noncommutative
residue of a holomorphic dif\/ferential operator $\mathcal D$ on $\cale$ for the case when $\cale$ is an arbitrary vector bundle over an
arbitrary compact complex manifold $X$. Theorem~\ref{th3} easily implies that this number is $C$ times the Lefschetz number of $\cale$ ($C$ being the
constant from Theorem~\ref{th3}). We also remark that if we replace $\Psi_{2n+1}$ by a $GL(n)$-basic cocycle $\Theta_{2n+1}$ representing the cohomology class $[\Psi_{2n+1}]$, Theorem~\ref{th3} holds in full generality (i.e., for arbitrary vector bundles with bounded geometry over complex manifolds that are not parallelizable). Indeed, the proof of Lemma~5.2 of~[10] can be modif\/ied to show the existence of such a $GL(n)$-basic representative of $[\Psi_{2n+1}]$. This settles Conjecture~3.3 of \cite{Sh1}.

\textbf{Convention 0.} We now state some conventions that stand throughout this paper.
\begin{enumerate}\itemsep=0pt
\item[(i)] The term {\it complex} shall be understood to mean a complex with homological grading with one exception that we point out in (iv) below.

\item[(ii)] The term {\it double complex} will mean ``complex of complexes''. Let $\{C_{p,q}\}$ be a double
complex with horizontal and vertical dif\/ferentials $d_h$ and $d_v$ respectively. One associates two total complexes $\text{Tot}^{\oplus}(C)$ and
$\text{Tot}^{\Pi}(C)$ with the double complex $C$, where
\begin{gather*}
\big(\text{Tot}^{\oplus}(C)\big)_n:=\oplus_{p+q=n} C_{p,q}, \qquad \big(\text{Tot}^{\Pi}(C)\big)_n:=\Pi_{p+q=n} C_{p,q}.
\end{gather*}
The dif\/ferential we use on the shifted total complexes
 $\text{Tot}^{\oplus}(C)[k]$ and  $\text{Tot}^{\Pi}(C)[k]$ will satisfy $d(x)=d_h(x)+{(-1)}^{|x|}d_v(x)$ for any homogenous element $x$ of $\text{Tot}(C)[k]$.

\item[(iii)] The DG-algebras that we work with naturally come with cohomological grading. If $A^{\bullet}$ is a~(cohomologically graded) DG-algebra, we automatically view it as a homologically graded algebra by inverting degrees (i.e., elements of $A^n$ have chain degree $-n$). The notation $A^{\bullet}$ in this case is still continued in order not to introduce unnatural notation. The same principle holds for DG-modules over DG-algebras that we work with.

\item[(iv)] When computing the hypercohomology of a complex $\calg F_{\bullet}$ of sheaves on $X$, we convert $\calg F_{\bullet}$ into a complex with cohomological grading by inverting degrees.

\item[(v)] For a complex $\calg F_{\bullet}$ of sheaves on $X$, $\mathcal{H}_{\bullet}(\calg F)$ shall denote the (graded) sheaf of homologies of $\calg F_{\bullet}$.
\end{enumerate}

\section[The Hochschild chain complex $\Gamma(X,\widetilde{\text{hoch}(\dif(\cale))})$ and other preliminaries]{The Hochschild chain complex $\boldsymbol{\Gamma(X,\widetilde{\text{hoch}(\dif(\cale))})}$\\ and other preliminaries}\label{section2}

\textbf{2.1.} Let $\difdt(\cale)(U)$ denote $\Gamma(U,\Omega^{0,\bullet}_X \otimes_{\strcc{X}} \dif(\cale))$. Let $A^{\bullet}$ be any
DG-algebra. Let $\text{C}_{\bullet}(A^{\bullet})$ denote the complex of Hochschild chains of $A^{\bullet}$. Whenever necessary,
$\text{C}_{\bullet}(A^{\bullet})$ will be converted to a cochain complex by inverting degrees. Note that the dif\/ferential on
$\text{C}_{\bullet}(\difdt(\cale)(U))$ extends to a dif\/ferential on the graded vector space $\oplus_k \difdt(\cale^{\boxtimes k})(U^k)[k-1]$ where $\cale^{\boxtimes k}$ denotes the vector bundle $\pi_1^* \cale \otimes \cdots \otimes \pi_k^*\cale$ on $X^k$ (with $\pi_j:X^k \rar X$ denoting the projection to the $j$-th factor).
We denote the resulting complex by $\widehat{\text{C}_{\bullet}}(\difdt(\cale)(U))$.

Let $\widetilde{\text{hoch}(\dif(\cale))}$ denote the sheaf associated to the presheaf
\begin{gather*}
U \mapsto
\widehat{\text{C}_{\bullet}}(\difdt(\cale)(U))
 \end{gather*}
 of complexes of $\compl$-vector spaces. Note that left multiplication of a homogenous element of  $\difdt(\cale)(U)$ by a smooth function on $U$ is well def\/ined. It follows that each $\difdt(\cale^{\boxtimes k})(U^k)$ is a~graded left $\text{C}^{\infty}(U)$-module via
\begin{gather*}
f.D := \pi_1^*(f).D \quad \text{for all} \ \ f \in \text{C}^{\infty}(U),\qquad  D \in \difdt(\cale^{\boxtimes k})(U^k).
\end{gather*}
Hence,  $\widetilde{\text{hoch}(\dif(\cale))}$ is a
complex each of whose terms is a module over the sheaf of smooth functions on $X$ (though the dif\/ferential is not compatible with this $C^{\infty}$-module structure). It follows that $\widetilde{\text{hoch}(\dif(\cale))}$ is a
complex of f\/labby sheaves on $X$. Therefore,
\begin{gather} \label{hypercoh1}
\text{H}_{\bullet}(\Gamma(X,\widetilde{\text{\rm hoch}(\dif(\cale))})) \simeq {\mathbb H}^{\bullet}(X,\widetilde{\text{hoch}(\dif(\cale))}).
\end{gather} {\it Here, one must convert the chain complex $\widetilde{\text{\rm hoch}(\dif(\cale))}$ into a cochain complex by
inverting degrees before taking hypercohomology.}

\begin{proposition} \label{prelim*}
\begin{gather*}
\text{\rm H}_{i}(\Gamma(X,\widetilde{\text{\rm hoch}(\dif(\cale))})) \simeq \text{\rm H}^{2n-i}(X,\compl)\qquad  \text{\rm whenever} \ \
0 \leq i \leq 2n.
\end{gather*}
 \end{proposition}

\begin{proof} Note that the dif\/ferential on the Hochschild chain complex $\text{C}_{\bullet}(\Gamma(U,\dif(\cale)))$ extends to a dif\/ferential
on the graded vector space $\oplus_k \Gamma(U^k,\dif(\cale^{\boxtimes k}))[k-1]$. Call the resulting complex
$\widehat{\text{C}_{\bullet}}(\Gamma(U,\dif(\cale)))$. One can then consider the sheaf $\widehat{\text{C}_{\bullet}}(\dif(\cale))$ associated
to the presheaf
\begin{gather*}
U \mapsto \widehat{\text{C}_{\bullet}}(\Gamma(U,\dif(\cale)))
\end{gather*}
 of complexes of $\compl$-vector spaces. We claim that the
natural map of complexes from $\widehat{\text{C}_{\bullet}}(\dif(\cale))$ to $\widetilde{\text{hoch}(\dif(\cale))}$ is a quasiisomorphism. To prove our claim, it suf\/f\/ices to verify that the natural map of complexes from $\widehat{\text{C}_{\bullet}}(\dif(\cale))$ to $\widetilde{\text{hoch}(\dif(\cale))}$ locally induces an isomorphism on cohomology sheaves. Indeed, on any open subset $U$ of $X$ on which $\cale$ and $T^{1,0}_X$ are trivial, the natural map from $\Gamma(U^k,\dif(\cale^{\boxtimes k}))$ to
$\Gamma(U^k,\difdt(\cale^{\boxtimes k}))$ is a quasiisomorphism (since $\strcc{X^k}$ is quasiisomorphic to $(\Omega^{0,\bullet}_{X^k},\bar{\partial})$ and $\Gamma(U^k,\dif(\cale^{\boxtimes k}))$ is a free module over~$\Gamma(U^k,\strcc{X^k})$). The desired proposition then follows from~\eqref{hypercoh1} and the fact that
$\widehat{\text{C}_{\bullet}}(\dif(\cale))$ is quasiisomorphic to the shifted constant sheaf~$\underline{\compl}[2n]$. This fact (see
\cite[Lemma~3]{Ram1}) is also implicit in many earlier papers, for instan\-ce,~\cite{Bryl}.
\end{proof}

\textbf{Convention 1.} {\it Whenever we identify $\text{\rm H}_{i}(\Gamma(X,\widetilde{\text{\rm hoch}(\dif(\cale))}))$ with
$\text{\rm H}^{2n-i}(X,\compl)$, we use the identification coming from a specific quasiisomorphism between
$\widehat{\text{\rm C}_{\bullet}}(\dif(\cale))$ and $\underline{\compl}[2n]$.} This quasiisomorphism is constructed as follows. ``Morita
invariance'' yields a quasiisomorphism between $\widehat{\text{C}_{\bullet}}(\dif(\cale))$ and $\widehat{\text{C}_{\bullet}}(\dif(X))$ (see the
proof of Lemma~3 of~\cite{Ram1}). By a result of~\cite{Bryl}, $\widehat{\text{C}_{\bullet}}(\dif(X))$ is quasiisomorphic to
$\underline{\compl}[2n]$. The quasiisomorphism we shall use throughout this paper is the one taking the class of the normalized Hochschild
$2n$-cycle
\begin{gather*}
\sum_{\sigma \in S_{2n}} \text{sgn}(\sigma)\left(1,z_1,\frac{\partial}{\partial z_1},\dots,z_n,\frac{\partial}{\partial z_{n}}\right)
\end{gather*} to
$1[2n]$ on any open subset $U$ of $X$ with local holomorphic coordinates $z_1,\dots,z_n$. We remark that in the above formula, any permutation
$\sigma$ in $S_{2n}$ permutes the last $2n$ factors leaving the f\/irst one f\/ixed.
 {\it We denote this identification
by $\beta_{\cale}$ as we did in}~\cite{Ram1}.

\textbf{2.1.1.} As in \cite{Ram1}, we def\/ine a topology on $\diffdt(\cale)$ as follows. Let $\text{Dif\/f}^{\leq k,\bullet}(\cale)$ denote
$\Gamma(X,\Omega^{0,\bullet}_X \otimes \dif^{\leq k}(\cale))$ where $\dif^{\leq k}(\cale)$ denotes the sheaf of holomorphic dif\/ferential
operators on $\cale$ of order $\leq k$. Equip $\cale$ and $\Omega^{0,\bullet}_X$ with Hermitian metrics. Equip $\text{Dif\/f}^{\leq
k,\bullet}(\cale)$ with the topology generated by the family of seminorms $\{||\cdot ||_{\phi,K,s}  \ | \, K \subset X \text{ compact}, \ s \in
\Gamma(K,\cale \otimes \Omega^{\bullet})$, $\phi \text{ a C}^{\infty}  \ \text{dif\/ferential
 operator on } \cale \}$ given by
\begin{gather*}
||D||_{\phi,K,s} = \text{Sup}\{||\phi D(s)(x)|| \ | \, x \in K\}.
\end{gather*}
 The topology on $\diffdt(\cale)$ is the direct limit of the
topologies on the $\text{Dif\/f}^{\leq k, \bullet}(\cale)$. More generally, for any open subset $U$ of $X$, one can def\/ine a topology on
$\Gamma(U,\difdt(\cale))$ in the same way.

\textbf{2.2.} Let $\cale$ be a holomorphic vector bundle on an {\it arbitrary} connected complex manifold $X$. Let $\Delta_{\cale}$ denote the
Laplacian of $\cale$ (for the operator $\bar{\partial}_{\cale}$). This depends on the choice of Hermitian metric for $\cale$ as well as on a
choice of a Hermitian metric for $X$. Recall that $\Delta_{\cale}=\Delta^{\cale}+F$ where $\Delta^{\cale}$ is the Laplacian of a $\text{C}^{\infty}$-connection on
$\cale$ (see Def\/inition 2.4 of \cite{BGV}) and $F \in \Gamma(X,\text{End}(\cale))$.

\begin{definition}
%\textbf{Definition.}
We say that $\cale$ has {\it bounded geometry} if for some choice of Hermitian metric on~$\cale$, there exists a
connection $\triangledown_{\cale}$ on~$\cale$ such that $\Delta_{\cale}=\Delta^{\cale}+F$, where $\Delta^{\cale}$ is the Laplacian of~$\triangledown_{\cale}$ and $F \in \Gamma(X,\text{End}(\cale))$, and all covariant derivatives of the curvature of~$\triangledown_{\cale}$ as
well as those of~$F$ are bounded on~$X$. In particular, any holomorphic vector bundle on a compact complex manifold has bounded geometry.
\end{definition}

Let $\cale$ be a vector bundle with bounded geometry on $X$. Let $\dol{\cale}{L^2}$ denote the Hilbert space of square integrable sections of
$\dol{\cale}{}\!\!$. Then, $\text{e}^{-t\Delta_{\cale}}$ makes sense as an integral operator on $\dol{\cale}{L^2}$ for any $t>0$ (see~\cite{Don}). Let $\Gamma_c$ denote the functor ``sections with compact support''. Suppose that $\alpha \in
\Gamma_c(X,\widetilde{\text{hoch}(\dif(\cale))})$ is a $0$-chain. Let $\alpha_0$ denote the component of $\alpha$ with antiholomorphic degree
$0$.  We recall the main results of \cite{Ram1} and \cite{Ram2} in the following theorem.

\begin{theorem} \label{oldrecall2} \qquad

$(1)$ For any $0$-chain $\alpha$, $\alpha_0\text{e}^{-t\Delta_{\cale}}$ is a trace class operator on
$\dol{\cale}{L^2}$ for any $t>0$.

$(2)$ Further, if $\alpha$ is a $0$-cycle,
\begin{gather*}
\lim_{t \rar \infty}
\text{\rm str}\big(\alpha_0\text{\rm e}^{-t\Delta_{\cale}}\big) = \frac{1}{(2 \pi i)^n}\int_X [\alpha],
\end{gather*}
 where $[\alpha]$ denotes the class of $\alpha$ in
$\text{\rm H}^{2n}_c(X,\compl)$.
\end{theorem}

\begin{remark}\label{remark1}
 Let us recall some aspects of \cite{Ram1} and \cite{Ram2}.
 When $X$ is compact, the linear functional
 \begin{gather} \label{fls} \alpha \mapsto \lim_{t \rar \infty}
\text{str}\big(\alpha_0\text{e}^{-t\Delta_{\cale}}\big)
\end{gather} on $\Gamma_c(X,\widetilde{\text{hoch}(\dif(\cale))})$ is really an extension of
the Feigin--Losev--Shoikhet--Hochschild $0$-cocycle of $\diffdt(\cale)$. We therefore, denote the linear functional given by~\eqref{fls} by
$\fls^{\text{hoch}}$. To show that the Feigin--Losev--Shoikhet cocycle indeed gives a linear functional on
$\text{H}_0(\Gamma(X,\widetilde{\text{hoch}(\dif(\cale))}))$ takes the sole computational ef\/fort in this whole program. The crux of this is a
certain ``estimate'' of the Feigin--Losev--Shoikhet cocycle (Proposition~3 of \cite{Ram1}). This estimate is proven in ``greater than usual''
detail in~\cite{Ram1}. It is done using nothing more than the simplest heat kernel estimates (Proposition~2.37 in~\cite{BGV}).
\end{remark}

 The second part of Theorem~\ref{oldrecall2} can be exploited to extend the Feigin--Losev--Shoikhet cocycle to a linear functional on certain homology theories that are
closely related to $\widetilde{\text{hoch}(\dif(\cale))}$.

\textbf{2.3. Cyclic homology.} One obtains a ``completed'' version of the Bar complex of\linebreak  $\Gamma(U,\difdt(\cale))$ by replacing
$\Gamma(U,\difdt(\cale))^{\otimes k}$ by $\Gamma(U^k,\difdt(\cale^{\boxtimes k}))$. Sheaf\/if\/ication then gives us a complex
$\widetilde{\text{bar}(\dif(\cale))}$ of sheaves on $X$ having the same underlying graded vector space as
$\widetilde{\text{hoch}(\dif(\cale))}$. One can then construct Tsygan's double complex
\begin{gather} \label{Tsyg}
\begin{CD}
\cdots @<{1-t}<< \cdots @<N<< \cdots @<{1-t}<< \cdots \\
@VVV      @VVV    @VVV  \\
\widetilde{\text{hoch}(\dif(\cale))}_2 @<1-t<< \widetilde{\text{bar}(\dif(\cale))}_2 @<N<< \widetilde{\text{hoch}(\dif(\cale))}_2 @<1-t<<\cdots\\
@VVV      @VVV    @VVV  \\
\widetilde{\text{hoch}(\dif(\cale))}_1 @<1-t<< \widetilde{\text{bar}(\dif(\cale))}_1 @<N<<
\widetilde{\text{hoch}(\dif(\cale))}_1 @<1-t<<\cdots\\ @VVV      @VVV    @VVV  \\ \widetilde{\text{hoch}(\dif(\cale))}_0 @<1-t<<
\widetilde{\text{bar}(\dif(\cale))}_ 0 @<N<< \widetilde{\text{hoch}(\dif(\cale))}_0 @<1-t<<\cdots\\
@VVV      @VVV    @VVV  \\
\cdots @<{1-t}<< \cdots @<N<< \cdots @<{1-t}<< \cdots \\
\end{CD}
\end{gather}

The total complex of this double complex is denoted by $\widetilde{\text{Cycl}(\dif(\cale))}$. Here, the operator $t$ acts on
$\Gamma(U,\difdt(\cale))^{\otimes k}$ by
\begin{gather*}
 D_1 \otimes \cdots \otimes D_k \mapsto {(-1)}^k {(-1)}^{d_k(d_1+\cdots +d_{k-1})}D_k \otimes D_1
\otimes\cdots \otimes D_{k-1},
\end{gather*}
 where the $D_i$ are homogenous elements of $\Gamma(U,\difdt(\cale))$ of degree $d_i$. $t$ then extends to
$\Gamma(U^k,\difdt(\cale^{\boxtimes k}))$ by continuity. It further extends to a map from $\widetilde{\text{bar}(\dif(\cale))}$ to
$\widetilde{\text{hoch}(\dif(\cale))}$ by sheaf\/if\/ication. Similarly, $N$ is the extension of the endomorphism $1+t+t^2+\cdots +t^{k-1}$ of
$\Gamma(U,\difdt(\cale))^{\otimes k}$ to a endomorphism of the graded vector space $\widetilde{\text{hoch}(\dif(\cale))}$. One can also
consider the Connes complex $\widetilde{\text{Co}(\dif(\cale))}:= \widetilde{\text{hoch}(\dif(\cale))}/(1-t)$. There is a natural surjection
$\Pi:\widetilde{\text{Cycl}(\dif(\cale))} \rar \widetilde{\text{Co}(\dif(\cale))}$ given by the quotient map on the f\/irst column of the double
complex \eqref{Tsyg} and $0$ on other columns. We omit the proof of the following standard proposition.

\begin{proposition} \label{pqsm} $\Pi$ is a quasi-isomorphism.
 \end{proposition}

{\sloppy
Applying the linear functional $\fls^{\text{hoch}}$ on the f\/irst column of the double complex of compactly supported global sections of the
double complex \eqref{Tsyg}, one obtains a linear functional $\fls^{\text{Cycl}}$ on $\Gamma_c(X,\widetilde{\text{Cycl}(\dif(\cale))})$.
Moreover, the endomorphism $1-t$ of the vector space $\Gamma_c(X,\widetilde{\text{hoch}(\dif(\cale))}_0)$ kills the direct summand
$\diff^0(\cale)$. By def\/inition, $\fls^{\text{hoch}}(\alpha)$ depends only on the component of~$\alpha$ in $\diff^0(\cale)$. It follows that
$\fls^{\text{hoch}}$ vanishes on the image of $1-t$. It therefore, induces a~linear functional on
$\Gamma_c(X,\widetilde{\text{Co}(\dif(\cale))})$ which we will denote by $\fls^{\text{Co}}$. More explicitly, if~$\alpha$ is a~$0$-chain in
$\Gamma_c(X,\widetilde{\text{Co}(\dif(\cale))})$, $\fls^{\text{Co}}(\alpha)=\lim_{t \rar \infty}
\text{str}(\alpha_0\text{e}^{-t\Delta_{\cale}})$ where~$\alpha_0$ is the component of~$\alpha$ in~$\diff^0(\cale)$. The following proposition
is a consequence of part~2 of Theorem~\ref{oldrecall2}.

}

\begin{proposition}
$\fls^{\text{\rm Cycl}}$ and $\fls^{\text{\rm Co}}$ descend to linear functionals on
$\text{\rm H}_0(\Gamma_c(X,\widetilde{\text{\rm Cycl}(\dif(\cale))}))$ and $\text{\rm H}_0(\Gamma_c(X,\widetilde{\text{\rm Cycl}(\dif(\cale))}))$
respectively.
\end{proposition}

\begin{proof} We have already observed that $\fls^{\text{hoch}}$ vanishes on the image of $1-t$ in
\linebreak $\Gamma_c(X,\widetilde{\text{hoch}(\dif(\cale))}_0)$. That it also vanishes on the image of the Hochschild boundary $b$ of
$\Gamma_c(X,\widetilde{\text{hoch}(\dif(\cale))})$ is implied by part 2 of Theorem \ref{oldrecall2}. The desired proposition then follows from
the construction of $\fls^{\text{Cycl}}$ and $\fls^{\text{Co}}$. \end{proof}

We recall from \cite{Ram2} that
\begin{gather*}
\text{H}_{0}(\Gamma_c(X,\widetilde{\text{Cycl}(\dif(\cale))})) \simeq \text{H}^{2n}(X,\compl).
\end{gather*}
It follows from Proposition \ref{pqsm} that
\begin{gather*}
\text{H}_{0}(\Gamma_c(X,\widetilde{\text{Co}(\dif(\cale))})) \simeq \text{H}^{2n}(X,\compl)
\end{gather*}
as well. Denote the natural projection from $\Gamma_c(X,\widetilde{\text{hoch}(\dif(\cale))})$ to
$\Gamma_c(X,\widetilde{\text{Co}(\dif(\cale))})$ by $p$. The following proposition is immediate from the def\/inition of $\fls^{\text{Co}}$.

\begin{proposition} \label{nontriv1} The following diagram commutes.
\begin{gather*}
\begin{CD} \text{\rm H}_0(\Gamma_c(X,\widetilde{\text{\rm hoch}(\dif(\cale))}))
@>\text{\rm H}_0(p)>>\text{\rm H}_0(\Gamma_c(X,\widetilde{\text{\rm Co}(\dif(\cale))}))\\ @VV{\fls^{\text{\rm hoch}}}V  @V{\fls^{\text{\rm Co}}}VV\\ \compl  @>
\text{\rm id}>> \compl\\ \end{CD}
\end{gather*}
\end{proposition}

Note that if $\beta$ is a nontrivial $0$-cycle of $\Gamma_c(X,\widetilde{\text{hoch}(\dif(\cale))})$, the above proposition and Theo\-rem~\ref{oldrecall2} together imply that $\fls^{\text{Co}}(p(\beta))=\frac{1}{(2 \pi i)^n}\int_X [\beta] \neq 0$. Therefore, $\fls^{\text{Co}}$ is
a nontrivial linear functional on $\text{H}_0(\Gamma_c(X,\widetilde{\text{Co}(\dif(\cale))}))$.

\textbf{2.4.~Lie algebra homology.} We begin this subsection by reminding the reader of a convention that we will follow.

\textbf{Convention~2.} Let $\mathfrak{g}$ be any DG Lie algebra with dif\/ferential of degree 1. By the term ``complex of Lie chains of
$\mathfrak{g}$'', we shall refer to the total complex of the double complex $\{C_{p,q}:=(\wedge^p\mathfrak{g})_{-q} | p \geq 1 \}$. The
horizontal dif\/ferential is the Chevalley--Eilenberg dif\/ferential and the vertical dif\/ferential is that induced by the dif\/ferential on
$\mathfrak{g}$. Note that our coef\/f\/icients are in the trivial $\mathfrak{g}$-module, and that {\it we are brutally truncating the column of
the usual double complex that lies in homological degree $0$}. In particular, if $\mathfrak{g}$ is concentrated in degree $0$, we do not
consider Lie $0$-chains of $\mathfrak{g}$.

\textbf{2.4.1.} Recall that $\difdt(\cale)$ is a sheaf of topological algebras on $X$. Over each open subset $U$ of $X$, the topology on
$\Gamma(U,\difdt(\cale))$ is given as in Section~2.1.1. It follows that $\mathfrak{gl}_r(\difdt(\cale))$ is a sheaf of topological Lie
algebras on $X$ with the topology on $\Gamma(U,\mathfrak{gl}_r(\difdt(\cale)))$ induced by the topology on $\Gamma(U,\difdt(\cale))$ in the
obvious way. Treating $\Gamma(U,\fmat(\difdt(\cale)))$ as the direct limit of the topological vector spaces
$\Gamma(U,\mathfrak{gl}_r(\difdt(\cale)))$ makes $\fmat(\difdt(\cale))$ a~sheaf of topological Lie algebras on $X$. In this case, the complex
of ``completed'' Lie chains $\widehat{\text{C}}^{\rm lie}_{\bullet}(\Gamma(U,\fmat(\difdt(\cale))))$ is obtained by replacing
$\wedge^k\Gamma(U,\fmat(\difdt(\cale)))$ by the image of the idempotent $\frac{1}{k!}\sum_{\sigma \in S_k} \text{sgn}(\sigma)\sigma$ acting on
$\fmat(\compl)^{\otimes k} \otimes \Gamma(U^k,\difdt(\cale^{\boxtimes k}))$. Here, $\sigma$ acts on the factors $\fmat(\compl)^{\otimes k}$
and $\Gamma(U^k,\difdt(\cale^{\boxtimes k}))$ simultaneously. In an analogous fashion, one can construct the complex
$\widehat{\text{C}}^{\rm lie}_{\bullet}(\Gamma(U,\fmat(\dif(\cale))))$.
 Note that $\widehat{\text{C}}^{\rm lie}_{\bullet}(\Gamma(U,\fmat(\dif(\cale))))$ is quasi-isomorphic
 to $\widehat{\text{C}}^{\rm lie}_{\bullet}(\Gamma(U,\fmat(\difdt(\cale))))$ since the columns of the corresponding double complexes are
 quasi-isomorphic.

\begin{definition}
 We def\/ine $\widetilde{\text{Lie}}(\dif(\cale))$ to be the complex
of sheaves associated to the complex of presheaves
\begin{gather*}
 U \mapsto
  \widehat{\text{\rm C}}^{\rm lie}_{\bullet}(\Gamma(U,\fmat(\difdt(\cale))))
  \end{gather*}
  of $\compl$-vector spaces on $X$.
\end{definition}

  In the proposition below, the subscript `$c$' indicates ``compact supports''.

\begin{proposition}\label{proposition5}
\begin{gather*}
\text{\rm H}_i(\Gamma_c(X,\widetilde{\text{\rm Lie}}(\dif(\cale)))) \simeq \mathbb{H}^{-i}_c(\text{\rm Sym}^{\geq
1}(\underline{\compl}[2n+1] \oplus \underline{\compl}[2n+3] \oplus \cdots )).
\end{gather*}
\end{proposition}
\begin{proof} Recall from \cite{FT} as well as Section 10.2 of \cite{lod} that for any $\compl$-algebra $A$, we have an isomorphism
 \begin{gather*}
 \text{H}_{\geq 1}^{\rm lie}(\fmat(A)) \simeq \text{Sym}^{\geq 1}(\text{HC}_{\bullet}(A)[1]).
\end{gather*}
 Section 10.2 of \cite{lod} shows the following facts
 \begin{enumerate}\itemsep=0pt
\item[(i)] The complex $\text{C}^{\rm lie}_{\bullet}(\fmat(A))_{\fmat(\compl)}$ is a commutative and cocommutative DG-Hopf-algebra.

\item[(ii)] One has a quasiisomorphism from the primitive part $\text{Prim}(\text{C}^{\rm lie}_{\bullet}(\fmat(A))_{\fmat(\compl)})$ to the Connes complex $\text{C}^{\lambda}_{\bullet}(A)[1]$.

\item[(iii)] The natural map of complexes $\text{C}^{\rm lie}_{\bullet}(\fmat(A)) \rar \text{C}^{\rm lie}_{\bullet}(\fmat(A))_{\fmat(\compl)}$ is a quasi-isomorphism.
    \end{enumerate}

 It follows from these facts (and the Milnor--Moore and PBW theorems) that one has the quasi-isomorphisms
 \begin{gather*}
 \text{C}^{\rm lie}_{\bullet}(\fmat(A)) \rar \text{C}^{\rm lie}_{\bullet}(\fmat(A))_{\fmat(\compl)} \leftarrow \text{Sym}^{\geq 1}(\text{Prim}(\text{C}^{\rm lie}_{\bullet}(\fmat(A))_{\fmat(\compl)})),\\
 \text{Sym}^{\geq 1}(\text{Prim}(\text{C}^{\rm lie}_{\bullet}(\fmat(A))_{\fmat(\compl)})) \rar \text{Sym}^{\geq 1}(\text{C}^{\lambda}_{\bullet}(A)[1]),\\
 \text{Sym}^{\geq 1}(\text{CC}_{\bullet}(A)[1]) \rar \text{Sym}^{\geq 1}(\text{C}^{\lambda}_{\bullet}(A)[1]).
 \end{gather*}

 The reader should recall that we follow Convention 2 when def\/ining the Lie chain complex. The last quasi-isomorphism above is induced by the standard quasi-isomorphism from $\text{CC}_{\bullet}(A)$ to~$\text{C}^{\lambda}_{\bullet}(A)$.
 The discussion of Section~10.2 of \cite{lod} goes through with the obvious modif\/ications when the algebra $A$ is replaced the sheaf of topological algebras $\dif(\cale)$ and Lie algebra and cyclic chains are completed in the sense mentioned above (see Section~2.3.3 of \cite{bktv} for a similar assertion in the $\text{C}^{\infty}$ case). Recall from \cite{will} (see also \cite{wod,Ram1}) that the completed cyclic chain complex of $\dif(\cale)$ is quasi-isomorphic to $\underline{\compl}[2n] \oplus \underline{\compl}[2n+2] \oplus\cdots $. It follows that
$\widehat{\text{C}}^{\rm lie}_{\bullet}(\fmat(\dif(\cale)))$ is isomorphic to $\text{Sym}^{\geq 1}(\underline{\compl}[2n+1] \oplus
\underline{\compl}[2n+3] \oplus \cdots )$ in the derived category $\text{D}(\text{Sh}_{\compl}(X))$ of sheaves of $\compl$-vector spaces on $X$. Further note that
$\widehat{\text{C}}^{\rm lie}_{\bullet}(\fmat(\dif(\cale)))$ is quasi-isomorphic to
$\widetilde{\text{Lie}}(\dif(\cale))$ (since $\strcc{X}$ is quasi-isomorphic to $\Omega^{0,\bullet}_X$). Hence,
$\mathbb{H}^{-i}_c(\widetilde{\text{Lie}}(\dif(\cale)))$ is isomorphic to $\mathbb{H}^{-i}_c(\text{Sym}^{\geq 1}(\underline{\compl}[2n+1] \oplus
\underline{\compl}[2n+3] \oplus \cdots))$ for all $i$.

As in Section~2.1, one argues that $\widetilde{\text{Lie}}(\dif(\cale))$ is a complex of f\/labby sheaves on $X$. Hence, $\text{H}_i(\Gamma_c(X,\widetilde{\text{Lie}}(\dif(\cale))))$ is isomorphic to $\mathbb{H}^{-i}_c(\widetilde{\text{Lie}}(\dif(\cale)))$, proving the desired proposition.
 \end{proof}

Note that $\text{Sym}^{\geq 1}(\underline{\compl}[2n+1] \oplus \underline{\compl}[2n+3] \oplus \cdots)$ is a direct sum of copies of the
constant sheaf~$\underline{\compl}$ concentrated in homological degrees $\geq 2n+1$. Let $a(n,i)$ denote the number of partitions of~$i$ into pairwise distinct summands from the set $\{2n+1,2n+3,2n+5,\dots\}$. By Proposition 3.2.2 of \cite{KS},

\begin{corollary}\label{cor1}
$\text{\rm H}_i(\Gamma_c(X,\widetilde{\text{\rm Lie}}(\dif(\cale)))) \simeq \oplus_{j \in \mathbb Z} \text{\rm H}^{j}_c(X,\compl)^{a(n,i+j)}$ for all $i \in \mathbb Z$. In particular, $\text{\rm H}_1(\Gamma_c(X,\widetilde{\text{\rm Lie}}(\dif(\cale)))) \simeq \text{\rm H}^{2n}_c(X,\compl)$ and
$\text{\rm H}_i(\Gamma_c(X,\widetilde{\text{\rm Lie}}(\dif(\cale)))) \simeq 0$ $\forall\,  i \leq 0$.
 \end{corollary}

\textbf{2.4.2.} For any algebra $A$, we have a map of complexes
\begin{gather*}
{L}: \ \ \text{C}^{\rm lie}_{\bullet}(\fmat(A))[1] \rar
\text{C}^{\lambda}_{\bullet}(A),
\end{gather*} where the right hand side is the Connes complex of $A$ (the quotient of the Hochschild chain complex of $A$
by the image of $1-t$). To be explicit, a Lie $k+1$ chain $M_0 \wedge \cdots \wedge M_k$ maps to the $k$-chain $ \sum_{\sigma \in S_{k}} M_0
\otimes M_{\sigma(1)} \otimes \cdots \otimes M_{\sigma(k)}$ of $\text{C}^{\lambda}_{\bullet}(M_{N}(A))$ for $N$ suf\/f\/iciently large. One then
applies the generalized trace map to obtain a $k$ chain in $\text{C}^{\lambda}_{\bullet}(A)$.

{\sloppy
It is easily verif\/ied that the map $L$ extends by continuity to a map of complexes from
$\widehat{\text{C}}^{\rm lie}_{\bullet}(\Gamma(U,\fmat(\difdt(\cale))))[1]$ to the completed Connes complex for $\Gamma(U,\difdt(\cale))$. By
sheaf\/if\/ication followed by taking global sections, it further yields a map of complexes from
$\Gamma_c(X,\widetilde{\text{Lie}}(\dif(\cale))[1])$ to $\Gamma_c(X,\widetilde{\text{Co}(\dif(\cale))})$ which we shall continue to denote by~$L$.

}

\begin{proposition} \label{nontriv2} $L$ induces a nonzero map from
\begin{gather*}
\text{\rm H}_0(\Gamma_c(X,\widetilde{\text{\rm Lie}}(\dif(\cale))[1])) \qquad {\rm to}\qquad
 \text{\rm H}_0(\Gamma_c(X,\widetilde{\text{\rm Co}(\dif(\cale))})).
 \end{gather*} \end{proposition}

\begin{proof}\sloppy Recall from \cite{lod} that for any algebra $A$, $L$ induces an isomorphism between the primitive part of
$\text{H}^{\rm lie}_{\bullet}(\fmat(A))$ and $\text{HC}_{\bullet-1}(A)$. For the topological algebra $\Gamma(U,\difdt(\cale))$ on a~suf\/f\/iciently
small open subset $U$ of $X$, it induces an isomorphism between the primitive part of
$\text{H}_{\bullet}(\widehat{\text{C}}^{\rm lie}_{\bullet}(\Gamma(U,\fmat(\difdt(\cale)))))$ and the completed cyclic homology of
$\Gamma(U,\difdt(\cale))$. In particular, it induces an isomorphism between
\begin{gather*}
\text{H}_{2n}(\widehat{\text{C}}^{\rm lie}_{\bullet}(\Gamma(U,\fmat(\difdt(\cale))))[1]) \qquad {\rm and} \qquad
\text{H}_{2n}(\widehat{\text{C}}^{\lambda}_{\bullet}(\Gamma(U,\fmat(\difdt(\cale)))))
\end{gather*}
 both of which are isomorphic to $\compl$. The desired
proposition follows from this. \end{proof}

 It follows from
Propositions \ref{nontriv1} and \ref{nontriv2} that the map $\fls^{\text{Co}} \circ L$ induces a nonzero linear functional on
$\text{H}_0(\Gamma_c(X,\widetilde{\text{Lie}}(\dif(\cale))[1]))$. We will denote this linear functional by $\fls^{\rm lie}$ at the level of chains
as well as on homology. Note that if $\alpha$ is a $0$-chain in $\Gamma_c(X,\widetilde{\text{Lie}}(\dif(\cale))[1])$, the component of~$L(\alpha)$ in $\diff^0(\cale)$ is $\text{Tr}(\alpha_0)$ where $\alpha_0$ is the component of $\alpha$ in $\fmat(\diff^0(\cale))$ and
$\text{Tr}$ is the usual matrix trace. Therefore,
\begin{gather*} %\label{flslie}
\fls^{\rm lie}(\alpha)= \lim_{t \rar \infty}
\text{str}\big(\text{Tr}(\alpha_0)e^{-t\Delta_{\cale}}\big). \end{gather*}

Let $\calg F$ be another holomorphic vector bundle on $X$. There is a natural map $\iota:\difdt(\cale) \rar \difdt(\cale \oplus \calg F)$ of
DG-algebras on $X$. $\iota$ induces a map
\begin{gather*}
\bar{\iota}: \ \ \Gamma_c(X,\widetilde{\text{Lie}}(\dif(\cale))[1]) \rar
\Gamma_c(X,\widetilde{\text{Lie}}(\dif(\cale \oplus \calg F))[1]).
\end{gather*}

\begin{proposition} \label{book1} The following diagram commutes
\begin{gather*}
\begin{CD}
   \text{\rm H}_0(\Gamma_c(X,\widetilde{\text{\rm Lie}}(\dif(\cale))[1]))
   @>\text{\rm H}_0(\bar{\iota})>> \text{\rm H}_0(\Gamma_c(X,\widetilde{\text{\rm Lie}}(\dif(\cale \oplus \calg
   F))[1]))\\
   @VV{\fls^{{\rm lie},\cale}}V   @V{\fls^{{\rm lie},\cale \oplus \calg F}}VV\\
   \compl @>\text{\rm id}>> \compl\\
   \end{CD}
   \end{gather*}
\end{proposition}

\begin{proof} Let $\alpha$ be a $0$-chain in $\Gamma_c(X,\widetilde{\text{Lie}}(\dif(\cale))[1])$. Clearly,
$\text{Tr}(\bar{\iota}(\alpha)_0)= \iota(\text{Tr}(\alpha_0))$. It then follows immediately that
\begin{gather*}
 \lim_{t \rar \infty}
\text{str}\big(\text{Tr}(\bar{\iota}(\alpha_0))e^{-t\Delta_{\cale \oplus \calg F}}\big) = \lim_{t \rar \infty}
\text{str}\big(\text{Tr}(\alpha_0)e^{-t\Delta_{\cale}}\big)
\end{gather*}
 provided the Hermitian metric on $\cale \oplus \calg F$ is chosen by retaining that on
$\cale$ and choosing an arbitrary Hermitian metric on $\calg F$. Part~2 of Theorem {\ref{oldrecall2}} however, implies that on homology,
$\fls^{{\rm lie},\cale \oplus \calg F}$ is independent of the choice of Hermitian metric on $\cale \oplus \calg F$. This proves the desired
proposition.
\end{proof}

\section{Dif\/ferential forms computing the Lefschetz number}\label{section3}

\textbf{3.1. Fedosov dif\/ferentials.} The material in this subsection is a by and large a rehash of material from \cite{enfe} and \cite{Sh1}.
We extend a Gelfand-Fuks type construction from \cite{enfe} in this section. This construction is a version of a construction that was
originally done in \cite{GF} by I.M.~Gel'fand and D.B.~Fuks for Lie cocycles with trivial coef\/f\/icients for the Lie algebra of smooth vector
f\/ields on a smooth manifold. We begin by recalling the construction from~\cite{enfe}.

\textbf{3.1.1.} Let $r$ be the rank of $\cale$. Let $J_p\cale$ denote the bundle of $p$-jets of local trivializations of $\cale$. In
particular, $J_1\cale$ is the extended frame bundle whose f\/iber over each $x \in X$ consists of the set of all pairs comprising a basis of
$T^{1,0}_{X,x}$ and a basis of $\cale_x$. The group $G:=GL(n,\compl) \times GL(r,\compl)$ acts on the right on $J_p\cale$ for each $p$. More
specif\/ically, given a local isomorphism of bundles $\compl^n \times \compl^r \rar \cale$, $GL(n,\compl)$ acts by linear coordinate
transformations on $\compl^n$ and $GL(r,\compl)$ acts by linear transformations on $\compl^r$. This makes $J_1\cale$ a principal $G$-bundle.
Let $J_{\infty}\cale$ denote $\varprojlim J_p\cale$. Since the natural submersions $J_{p+1}\cale \rar J_p\cale$ are $G$-equivariant, the
natural map $J_{\infty}\cale \rar J_1\cale$ descends to a map $J_{\infty}\cale/G \rar J_1\cale/G=X$. Since the f\/ibres of this map are
contractible (see~\cite{enfe}), there exists a smooth section $\varphi$ of $J_{\infty}\cale/G$ over $X$. This is equivalent to a $G$-equivariant section of
$J_{\infty}\cale$ over $J_1\cale$. This section is unique upto smooth homotopy.

\textbf{3.1.2.} Let $\strcc{n}$ denote $\compl[[y_1,\dots,y_n]]$. Let $W_n$ denote the Lie algebra of formal holomorphic vector f\/ields (i.e.,
those of the form $\sum\limits_{i=1}^n f_i(y_1,\dots,y_n)\partial_i$ with $f_i \in \strcc{n}$ and $\partial_i:=\frac{\partial}{\partial y_i}$). Let
$W_{n,r}$ denote the Lie algebra $W_n \ltimes \mathfrak{gl}_r(\strcc{n})$. This semidirect product comes from the action of~$W_n$ on~$\mathfrak{gl}(\strcc{n})$ by derivations. We recall from~\cite{enfe} that there is a map of Lie algebras from $W_{n,r}$ to the Lie algebra of holomorphic vector f\/ields on~$J_{\infty}\cale$ which yields an isomorphism
$W_{n,r} \rar T^{1,0}_{\phi}J_{\infty}\cale$ for each $\phi \in J_{\infty}\cale$. This is equivalent to a~$W_{n,r}$-valued holomorphic
$1$-form $\Omega$ on $J_{\infty}\cale$ satisfying the Maurer--Cartan equation
\begin{gather*}
d\Omega+\frac{1}{2}[\Omega,\Omega]=0.
\end{gather*}
 Since there is a
natural map of Lie algebras $W_{n,r} \rar \text{M}_r(\diff_n)$, $\Omega$ may be viewed as a Maurer--Cartan form with values in
$\text{M}_r(\diff_n)$.

We clarify that for each $\phi=(\phi_p)$ a point in $J_{\infty}\cale$, $T^{1,0}_{\phi}J_{\infty}\cale:=\varprojlim T^{1,0}_{\phi_p}J_p\cale$ and $\Omega^{1,0}_{\phi} J_{\infty}\cale:=\varinjlim
\Omega^{1,0}_{\phi_p}J_p\cale$, etc. This is compatible with the def\/inition of $J_{\infty}\cale$ as $\varprojlim J_p\cale$. This enables us to make sense of tangent and cotangent spaces of the ``inf\/inite dimensional complex manifold''~$J_{\infty}\cale$, and hence, make sense of vector f\/ields and dif\/ferential forms on $J_{\infty}\cale$.

\textbf{3.1.3.} Consider the bundle $\dga{\cale}:= J_1\cale \times_G \text{M}_r(\diff_n)$ of algebras with f\/iberwise product. Any
trivialization of $J_1\cale$ over an open subset $U$ of $X$ yields an isomorphism of algebras
\begin{gather} \label{prelim}
\text{C}^{\infty}(U,\dga{\cale}) \simeq \text{C}^{\infty}(U,\text{M}_r(\diff_n)).
\end{gather}
More generally, a trivialization of
$J_1\cale$ over $U$ yields an isomorphism of graded algebras
\begin{gather*}
\Omega^{\bullet}(U,\dga{\cale})) \rar \Omega^{\bullet}(U,\text{M}_r(\diff_n)).
\end{gather*}
 The $G$ equivariant section $\varphi:J_1\cale \rar J_{\infty}\cale$ allows us to pull back $\Omega$. This gives us a one form
$\omega:=\varphi^*\Omega \in \Omega^1(J_1\cale,\text{M}_r(\diff_n))^G$. This in turn yields a f\/lat connection $D$ on the bundle of algebras~$\dga{\cale}$. Let us be explicit here. Given a trivialization of $J_1\cale$ over $U$, $\omega$ descends to a~$\text{M}_r(\diff_n)$ valued one
form $\omega_U$ on $U$ satisfying the Maurer--Cartan equation. The isomorphism~\eqref{prelim} identif\/ied the connection $D$ with the
``twisted'' de~Rham dif\/ferential $d+[\omega_U,-]$. More generally, the DG-algebra $(\Omega^{\bullet}(U,\dga{\cale}),D)$ is identif\/ied with the
DG-algebra $(\Omega^{\bullet}(U,\text{M}_r(\diff_n)),d+[\omega_U,-])$. Over any open subset $U$ of $X$, the degree $0$ sections $\hat{\mathcal
D}$ of $\Omega^{\bullet}(\dga{\cale})$ satisfying
\begin{gather*}
D{\hat{\mathcal D}}=0
\end{gather*} are in one-to-one correspondence with holomorphic dif\/ferential
operators on $\cale$ over $U$. The dif\/ferential $D$ is a {\it Fedosov differential} on the sheaf $\Omega^{\bullet}(\dga{\cale})$ of
DG-algebras on $X$.

\textbf{3.1.4.} Fedosov dif\/ferentials have earlier been constructed in a holomorphic setting in \cite{NT1}, followed by \cite{CDH} and more recently, \cite{CR}.
Construction of a Fedosov dif\/ferential $D$ on $\Omega^{\bullet}(\dga{\cale})$ has been explained in \cite{CR} (see also Dolgushev, \cite{dolg}) in the case when $\cale=\strcc{X}$. This goes through with minor modif\/ications for the case when
$\cale$ is arbitrary. This is a more ``careful'' construction of a~Fedosov dif\/ferential: it has the following special property (Part~2 of
Theorem~10.5 of~\cite{CR}: note that we do not need Part~1 of this result).

\begin{theorem} \label{fedosov} There is a map of sheaves of DG-algebras
\begin{gather*}
\difdt(\cale) \rar (\Omega^{\bullet}(\dga{\cale}),D).
\end{gather*}
\end{theorem}

Given that the construction of the Fedosov dif\/ferential $D$ from \cite{CR} is a priori dif\/ferent from that of \cite{enfe}, one would like to
know whether the Fedosov dif\/ferential $D$ from \cite{CR} can also be obtained via the construction from \cite{enfe} outlined earlier. The
reason why we need this relationship will become clear in later sections. The following paragraphs give a partial af\/f\/irmative answer to this
question that suf\/f\/ices for us. To be precise, the next paragraph shows that {\it there exists a $G$-equivariant section $\varphi:J_1\cale
\rar J_{\infty}\cale$ associated with any Fedosov differential on $\Omega^{\bullet}(\dga{\cale})$}. Further, Proposition~\ref{proposition8} in Section~3.1.5
shows that {\it for every $\mathcal{D} \in \Gamma(U, \dif(\cale))$, the section of~$\dga{\cale}$ over $U$ that $\mathcal D$ yields via
$\varphi$ $($via the construction in Section~$3.1.3)$ coincides with the flat $($with respect to $D)$ section of~$\dga{\cale}$ corresponding to
$\mathcal D$.} We also point out that on any open subset $U$ of $X$ on which $J_1\cale$ is trivial, any trivialization of $J_1\cale$ over $U$
identif\/ies the DG-algebra $(\Omega^{\bullet}(U,\dga{\cale}),D)$ with the DG-algebra $(\Omega^{\bullet}(U,\text{M}_r(\diff_n)),d+[\omega_U,-])$
for some $\text{M}_r(\diff_n)$-valued Maurer--Cartan form even in the construction of \cite{CR}.

\textbf{3.1.5.} The argument here follows Section~4.3 of \cite{fefesh}. Consider any point $(x,F)$ of $J_1\cale$ ($F$ being an extended frame
over $x$). For any holomorphic dif\/ferential operator ${\mathcal D}$ on $\cale$ over a~neighborhood of $x$, we have a unique section
$\hat{\mathcal D}$ of $\Omega({\mathcal B}_{\cale})$ satisfying $D{\hat{\mathcal D}}=0$. The frame $F$ identif\/ies $\hat{\mathcal D}$ with an
element of $\text{M}_r(\diff_n)$. This yields an isomorphism of algebras
\begin{gather*}
\Theta: \ \ \text{Jets}_x \dif(\cale) \rar \text{M}_r(\diff_n).
\end{gather*}
Moreover, since the above isomorphism preserves the order of the dif\/ferential operator, we obtain an isomorphism
\begin{gather*}
\theta:  \ \ \text{Jets}_x
\text{End}(\cale) \rar \text{M}_r({\mathcal O}_n).
\end{gather*}
 This corresponds to an element in the f\/ibre of $J_{\infty}\cale$ over $x$. We therefore,
 obtain a section $\varphi:J_1\cale \rar
J_{\infty}\cale$. It is easy to verify that this section is $G$-equivariant. Note that the isomorphism~$\theta$ induces an isomorphism
$\tilde{\theta}$ from $\text{Jets}_x \dif(\cale)$ to $\text{M}_r(\diff_n)$.

\begin{proposition}\label{proposition8}
 $\tilde{\theta}$ coincides with $\Theta$. \end{proposition}

 \begin{proof} Note that $\text{M}_r({\mathcal O}_n)$ and $\text{Id} \otimes
W_n$ generate $\text{M}_r(\diff_n)$ as an algebra. It therefore suf\/f\/ices to check that $\tilde{\theta}^{-1}$ and $\Theta^{-1}$ coincide on
$\text{M}_r({\mathcal O}_n)$ and $\text{Id} \otimes W_n$.Observe that the image of $W_{n,r}:=W_n \ltimes \mathfrak{gl}_r(\strcc{n})$ in
$\text{M}_r(\diff_n)$ is precisely the space of all elements of the form $X=D+f$ for some $D \in \text{Id} \otimes W_n$ and some $f \in
\text{M}_r({\mathcal O}_n)$. We shall continue to denote this image by $W_{n,r}$.

 \textbf{Claim.} Every derivation of $\text{M}_r(\strcc{n})$ is of the form
 $[X,-]$ where $X \in W_{n,r}$. Further, if $X \in
\text{M}_r(\diff_n)$ and $[X,-]$ is a derivation of $\text{M}_r(\strcc{n})$, then $X \in W_{n,r}$.

We f\/irst prove the proposition assuming our claim. Indeed, $\Theta^{-1}$ restricted to $\text{M}_r({\mathcal O}_n)$ coincides with
$\tilde{\theta}^{-1}$ restricted to $\text{M}_r({\mathcal O}_n)$ (both of which coincide with $\theta^{-1}$). For any $X \in W_{n,r}$ and for
any $g \in \text{M}_r(\strcc{n})$, $[\tilde{\theta}^{-1}(X),\theta^{-1}(g)]=\tilde{\theta}^{-1}([X,g])=\theta^{-1}([X,g])=
\Theta^{-1}([X,g])=
 [\Theta^{-1}X,\theta^{-1}g]$. It follows that for any $X \in
 W_{n,r}$, $\tilde{\theta}^{-1}(X)-\Theta^{-1}(X)$ commutes with all elements of $\text{Jets}_x \text{End}(\cale)$. In other words,
 $\tilde{\theta}^{-1}(X)-\Theta^{-1}(X)$ must be an element of $\theta^{-1}(\text{id} \otimes \strcc{n})$ for all
 $X \in W_{n,r}$. In particular, if  $X \in \text{id} \otimes W_n$, then $\tilde{\theta}^{-1}(X)=\Theta^{-1}(X)+\theta^{-1}(g)$ for some $g
 \in \text{id} \otimes \strcc{n}$.
 In this case,
 \begin{gather*}
 \tilde{\theta}^{-1}\big(X^2\big)= \Theta^{-1}\big(X^2\big)+ \Theta^{-1}(X).\theta^{-1}(g)+\theta^{-1}(g).\Theta^{-1}(X)+\theta^{-1}\big(g^2\big) \\
 \phantom{\tilde{\theta}^{-1}\big(X^2\big)}{}
 =\Theta^{-1}\big(X^2\big)+\Theta^{-1}(2gX)+\theta^{-1}\big([X,g]+g^2\big).
 \end{gather*}
 It follows that $\tilde{\theta}^{-1}(X^2)-\Theta^{-1}(X^2)$ is in $\theta^{-1}(\text{id} \otimes
 \strcc{n})$ only if $g=0$. This shows that $\tilde{\theta}^{-1}$
 and~$\Theta^{-1}$ coincide on $\text{id} \otimes W_n$ as well.

 Our claim remains to be proven. Recall (from \cite{lod} for instance) that for any algebra $A$,
\begin{gather*}
\text{HH}^1(A) \simeq \frac{\text{Der}(A,A)}{\text{Inner
 derivations}}.
 \end{gather*}
 In our case, the composite map
\begin{gather*}
W_n \rar \text{HH}^1(\strcc{n},\strcc{n}) \rar
 \text{HH}^1(\text{M}_r(\strcc{n}),\text{M}_r(\strcc{n}))
 \end{gather*} is an isomorphism.
 The f\/irst arrow above is the Hochschild--Kostant--Rosenberg map. The
 second arrow above is the cotrace. It follows that for any
 derivation $\phi$ of $\text{M}_r(\diff_n)$, there is an element $v$
 of $W_n$ such that the image of $v$ under the above composite map
 coincides with that of $\phi$ in
 $\text{HH}^1(\text{M}_r(\strcc{n}),\text{M}_r(\strcc{n}))$. Consider the
 derivation $\bar{v}:=[\text{id} \otimes v,-]$ of
 $\text{M}_r(\strcc{n})$. It is easy to verify that the image of
 $\phi-\bar{v}$ in $\text{HH}^1(\text{M}_r(\strcc{n}),\text{M}_r(\strcc{n}))$ is
 zero. Hence,  $\phi-\bar{v}$ is an inner derivation. In other
 words, as derivations on $\text{M}_r(\strcc{n})$, $\phi=[\text{id}
 \otimes v + f,-]$ for some $v \in W_n$ and some $f \in
 \text{M}_r(\strcc{n})$. This proves the f\/irst part of our claim.
 The  second part of our claim is then easy: if $[X,-]$ is a
 derivation of $\text{M}_r(\strcc{n})$, then $[X,-]=[Y,-]$ as derivations on $\text{M}_r(\strcc{n})$ for
 some
 $Y \in W_{n,r}$. Hence, $X-Y$ commutes with every element in
 $\text{M}_r(\strcc{n})$. This shows that $X-Y$ is an element of $\text{id} \otimes \strcc{n}$.
\end{proof}

\textbf{3.2. A Gel'fand--Fuks type construction.} We now recall and extend a construction from~\cite{enfe}.

\textbf{3.2.1.} Let $A$ be any DG-algebra. Let $\omega \in A^1$ be a Maurer--Cartan element (i.e., an satisfying $d\omega+\omega^2=0$). Let
$A_\omega$ denote the DGA whose underlying graded algebra is $A$ but whose dif\/ferential is $d+[\omega,-]$. Let $\text{C}_{\bullet}(A)$ denote
the complex of normalized Hochschild chains of $A$. Let $(\omega)^k$ denote the normalized Hochschild $0$-chain $(1,\omega,\dots,\omega)$
($k$ copies of $\omega$). Denote the set of $(p,q)$-shuf\/f\/les by $\text{Sh}_{p.q}$. There is a shuf\/f\/le product
\begin{gather*}
\times: \ \ A \otimes
{A}^{\otimes p} \otimes A \otimes A^{\otimes q} \rar A \otimes {A}^{\otimes p+q},\\
(a_0,a_1,\dots ,a_p) \otimes (b_0,a_{p+1},\dots,a_{p+q})\\
\qquad{} \mapsto \sum_{\sigma \in \text{Sh}_{p,q}} \text{sgn}(\sigma)(a_0b_0,a_{\sigma^{-1}(1)},\dots,a_{\sigma^{-1}(p)},a_{\sigma^{-1}(p+1)},\dots,
a_{\sigma^{-1}(p+q)}).
\end{gather*}

Note that $\text{C}_{\bullet}(A)= \oplus_{k \geq 1}A \otimes {A}^{\otimes k}[k-1]$ as graded vector spaces. The dif\/ferential on
$\text{C}_{\bullet}(A)$ extends naturally to a dif\/ferential on $\text{C}^{\Pi}_{\bullet}(A):= \prod\limits_{k \geq 1} A \otimes {A}^{\otimes k}[k-1]$
making the natural map of graded vector spaces a map of complexes. One has the following proposition due to Engeli and Felder \cite{enfe}.
Though \cite{enfe} used normalized Hochschild chains, their proof goes through in the current context as well.

\begin{proposition} \label{dgatwist} The map
\begin{gather*}
\text{\rm C}^{\Pi}_{\bullet}(A_{\omega}) \rar \text{\rm C}^{\Pi}_{\bullet}(A), \qquad x \mapsto x \times \sum_k
{(-1)}^k(\omega)^k
\end{gather*}
 is a term by term isomorphism of complexes. \end{proposition}

\textbf{3.2.2.} Recall that $GL(n,\compl)$ acts on a formal neighborhood of the origin by linear coordinate changes. This induces an action of
$GL(n,\compl)$ on $\diff_n$. The derivative of this action embeds the Lie algebra $\mathfrak{gl}_n$ in $\diff_n$ as the Lie subalgebra of
operators of the form $\sum_{j,k} a_{j,k}y_k\partial_j$. We say that the normalized cocycle $\alpha \in \text{C}^{p}(\diff_n)$ is {\it
$GL(n,\compl)$-basic} if $\alpha$ is $GL(n,\compl)$-invariant and if $\sum\limits_{i=1}^p {(-1)}^{i+1}\alpha(a_0,\dots,a_{i-1},a,a_{i+1},\dots,a_p)=0$ for
any $a \in \mathfrak{gl}_n$ and any $a_0,\dots,a_p \in \diff_n$.
 More generally, $G:=GL(n,\compl) \times GL(r,\compl)$ acts on $\text{M}_r(\diff_n)$. The action of $GL(n,\compl)$ on $\text{M}_r(\diff_n)$
 is induced by the action of $GL(n,\compl)$ on $\diff_n$ and $GL(r,\compl)$ acts on $\text{M}_r(\diff_n)$ by conjugation.
 The derivative of this action embeds $\mathfrak{g}:= \mathfrak{gl}_n \oplus \mathfrak{gl}_r$ as a Lie subalgebra of $\text{M}_r(\diff_n)$ of
 elements of the
form $\sum_{j,k}a_{j,k}z_k\partial_j \otimes \text{id}_{r \times r}+ B$ where
 $B \in \text{M}_r(\compl)$. Again, a cocycle $\alpha$ in $\text{C}^{p}(\text{M}_r(\diff_n))$ is said to be
 {\it $G$ basic} if $\alpha$ is $G$-invariant
 and $\sum\limits_{i=1}^p {(-1)}^{i+1}\alpha(a_0,\dots,a_{i-1},a,a_{i+1},\dots,a_p)=0$ for any $a \in \mathfrak{g}$ and any $a_0,\dots,a_p \in
 \text{M}_r(\diff_n)$. The following proposition is immediate from equation~\eqref{eq1}.

\begin{proposition} If $\psi_{2n}$ is $GL(n,\compl)$-basic then $\psi^r_{2n}$ is $G$-basic. \end{proposition}

  We also recall that a cocycle $\alpha \in \text{C}^{p}(\text{M}_r(\diff_n))$ is said to
  be {\it continuous} if it depends only on f\/initely many Taylor
  coef\/f\/icients of its arguments.

\textbf{3.2.3.} Let $\psi^r_{2n}$ be a continuous, $G$-basic Hochschild $2n$-cocycle of $\text{M}_r(\diff_n)$. Evaluation at~$\psi^r_{2n}$
gives a map of graded vector spaces
\begin{gather}
 \text{C}^{\Pi}_{\bullet}(\Omega^{\bullet}(U,M_r(\diff_n)),d) \rar \Omega^{2n-\bullet}(U),\nonumber\\
 (\omega_0a_0,\dots ,\omega_pa_p) \mapsto \psi^r_{2n}(a_0,\dots ,a_p)\omega_1\wedge\cdots \wedge \omega_p \nonumber\\
 \qquad  \text{for all} \ \ a_1,\dots ,a_p \in \text{M}_r(\diff_n), \ \ \omega_1,\dots ,\omega_p \in
\Omega^{\bullet}(U). \label{eval}
\end{gather}  By convention, $\psi^r_{2n}(a_0,\dots,a_p)=0$ if $p \neq 2n$.

\begin{proposition} \label{constr1} The map~\eqref{eval} is a map of complexes $($with $\Omega^{2n-p}(U)$ equipped with the diffe\-ren\-tial
$(-1)^pd_{\rm DR})$. \end{proposition}

\begin{proof} Note that any Hochschild $p$-chain $\alpha$ of $A:=(\Omega^{\bullet}(U,\text{M}_r(\diff_n)),d)$ can
be expressed as a~sum of terms $\sum_i \alpha_i$ where $\alpha_i$ has de~Rham degree $i$ and is an element of $A \otimes {A}^{p+i}$. By
def\/inition, $\psi^r_{2n}(\alpha)=\psi^r_{2n}(\alpha_{2n-p})$. Let $\delta$ denote the dif\/ferential on the Hochschild chain complex of $A$.
Then, $(\delta \alpha)_{2n-p}= {(-1)}^pd_{\rm DR}\alpha_{2n-p} + d_{\rm hoch}\alpha_{2n-p+1}$. Since $\psi^r_{2n}$ is a cocycle,
$\psi^r_{2n}(d_{\rm hoch}\alpha_{2n-p+1})=0$. The desired proposition is now immediate. \end{proof}

Observe that there is a natural map of complexes $ \text{C}_{\bullet}(\Omega^{\bullet}(U,{\mathcal B}_{\cale}),D) \rar
 \text{C}^{\Pi}_{\bullet}(\Omega^{\bullet}(U,{\mathcal B}_{\cale}),D)$.
 Also, the graded involution on $\Omega^{\bullet}(U)$ acting
 on $\Omega^{p}(U)$ by multiplication by ${(-1)}^{\lfloor \frac{p}{2} \rfloor}$
 is a map of complexes between $(\Omega^{2n-\bullet}(U),{(-1)}^{\bullet}d_{\rm DR})$ and
 $\Omega^{2n-\bullet}(U)$. Recall from Section~3.1.4 that a trivialization of $J_1\cale$ over $U$
 identif\/ies $(\Omega^{\bullet}(U,\dga{\cale}),D)$ with
 $\Omega^{\bullet}_{\omega}(U,\text{M}_{r}(\diff_n))$ where $\omega$
 is a $\text{M}_r(\diff_n)$-valued $1$-form on $U$ satisfying the
 Maurer--Cartan equation. By Propositions~\ref{dgatwist} and~\ref{constr1}, we have a map of complexes
\begin{gather} \label{B}
\text{C}_{\bullet}(\Omega^{\bullet}(U,{\mathcal B}_{\cale}),D) \rar \Omega^{2n-\bullet}(U).
\end{gather}

Explicitly, the map~\eqref{B} maps a chain $\mu \in \text{C}_{p}(\Omega^{\bullet}(U,{\mathcal B}_{\cale}),D)$ to
\begin{gather*}
{(-1)}^{\lfloor
\frac{2n-p}{2} \rfloor}\psi^r_{2n}\left(\hat{\mu} \times {\sum_k (-1)}^{k}(\omega)^{k}\right),
\end{gather*}
 where $\hat{\mu}$ is the image of $\mu$ in
$\text{C}^{\Pi}_p(\Omega^{\bullet}_{\omega}(U,\text{M}_r(\diff_n)))$.

\begin{proposition} \label{freetriv} The map \eqref{B} is independent of the trivialization of $J_1\cale$ used. \end{proposition}

\begin{proof} A
dif\/ferent trivialization of $J_1\cale$ over $U$ dif\/fers from the chosen one by a gauge change $g:U \rar G$. A section $\mu$ of
$\Omega^{\bullet}(U,{\mathcal B}_{\cale})$ transforms into $g.\mu$. The Maurer--Cartan form $\omega$ is replaced by $g.\omega-dg.g^{-1}$. That
$\psi^r_{2n}(\alpha \times \sum_k \pm (\omega)^k)=\psi^r_{2n}(g.\alpha \times \sum_k \pm (g.\omega-dg.g^{-1})^k)$ is immediate from the fact
that $\psi^r_{2n}$ is $G$-basic. This proves the desired proposition. \end{proof}

It follows that the map~\eqref{B} gives a map of complexes of presheaves
\begin{gather} \label{C} \big(U \mapsto
\text{C}_{\bullet}(\Omega^{\bullet}(U,{\mathcal B}_{\cale}),D)\big) \rar \big(U \mapsto \Omega^{2n-\bullet}(U)\big).
\end{gather} Note that
whatever we said so far is true for any Fedosov dif\/ferential no matter how it was constructed. Of course, the map of complexes of presheaves
\eqref{C} depends on the Fedosov dif\/fe\-ren\-tial~$D$. It follows from Theorem~\ref{fedosov} that we obtain a composite map of complexes of
presheaves
\begin{gather} \label{D} \big(U \mapsto \text{C}_{\bullet}(\difdt(\cale)(U))\big) \rar \big(U \mapsto
\text{C}_{\bullet}(\Omega^{\bullet}(U,{\mathcal B}_{\cale}),D)\big) \rar \big(U \mapsto \Omega^{2n-\bullet}(U)\big).
\end{gather}

\begin{proposition} \label{mapmade} The composite map \eqref{D} extends to a map of complexes of sheaves
\begin{gather*}
\widetilde{\text{\rm hoch}(\dif(\cale))}
\rar \Omega^{2n-\bullet}_X.
\end{gather*} \end{proposition}

\begin{proof} Denote the composite map \eqref{D} by $f_{\cale}$. For each open subset~$U$ of~$X$, $f_{\cale}$ yields a map of complexes from
$\text{C}_{\bullet}(\difdt(\cale)(U))$ to $\Omega^{2n-\bullet}(U)$. As described in Section~2.1.1,
$\widehat{\text{C}_{\bullet}}(\difdt(\cale)(U))$ has the structure of a (graded) locally convex topological vector space. One can verify
without dif\/f\/iculty that the subcomplex $\text{C}_{\bullet}(\difdt(\cale)(U))$ is dense in $\widehat{\text{C}_{\bullet}}(\difdt(\cale)(U))$ and
that the dif\/ferential of $\widehat{\text{C}_{\bullet}}(\difdt(\cale)(U))$ is continuous. Similarly, $\Omega^{2n-\bullet}(U)$ has the structure
of a (graded) complete locally convex topological vector space and the de~Rham dif\/ferential is continuous with respect to this topology.

Suppose that $f_{\cale}$ is continuous. Then, $f_{\cale}$ extends to a map of complexes from $\widehat{\text{C}_{\bullet}}(\difdt(\cale)(U)) $
to $\Omega^{2n-\bullet}(U)$ for each $U$. This can be easily see to yield a map of complexes of presheaves
\begin{gather*}
\big(U \mapsto
\widehat{\text{C}_{\bullet}}(\difdt(\cale)(U))\big) \rar \big(U \mapsto \Omega^{2n-\bullet}(U)\big).
\end{gather*}
 On sheaf\/if\/ication, $f$ yields a map of
complexes of sheaves
\begin{gather*}
\widetilde{\text{hoch}(\dif(\cale))} \rar \Omega^{2n-\bullet}_X
\end{gather*}
 proving the desired proposition. Indeed, continuity
of $f_{\cale}$ follows from continuity of $\psi^r_{2n}$. Since $\psi^r_{2n}$ is continuous, $C^l$-norms of $f(a_0 \otimes \cdots \otimes a_k)$
over a compact subsets of $U$ are estimated by $C^{l'}$-norms of f\/initely many Taylor coef\/f\/icients of the images $\hat{a_i}$ of the $a_i$ in
$\Omega^{\bullet}(U,\dga{\cale})$ over compact subsets of~$U$.
\end{proof}

{\it The map of complexes of sheaves obtained in Proposition~$\ref{mapmade}$ depends on the cocycle $\psi_{2n}$. We shall denote this map by
$f_{\cale,\psi_{2n}}$.} Recall that there is a natural map of complexes of sheaves $\widehat{\text{C}}_{\bullet}(\dif(\cale)) \rar
\widetilde{\text{hoch}(\dif(\cale))}$. Let $U$ be a subset of $X$ with local holomorphic coordinates $z_1,\dots,z_n$ such that $\cale$ is trivial
over $U$. By Convention~1, the standard generator for the homology of $\Gamma(U,\widehat{\text{C}}_{\bullet}(\dif(\cale)))$ is a cycle
$c_{\cale}(U)$ mapping to the normalized Hochschild $2n$-cycle
\begin{gather*}
\sum_{\sigma \in S_{2n}} \text{sgn}(\sigma)\sigma\left(\text{Id}_{\cale} \otimes
1, \text{Id}_{\cale} \otimes z_1, \text{Id}_{\cale} \otimes \frac{\partial}{\partial z_1},\dots, \text{Id}_{\cale} \otimes z_n,
\text{Id}_{\cale} \otimes \frac{\partial}{\partial z_n}\right).
\end{gather*}
 Let $c^r_{2n}$ denote the normalized Hochschild $2n$-cycle
 \begin{gather*}
 \sum_{\sigma
\in S_{2n}} \text{sgn}(\sigma)\sigma\left(\text{Id}_{r \times r}  \otimes 1, \text{Id}_{r \times r}  \otimes y_1, \text{Id}_{r \times r}  \otimes
\frac{\partial}{\partial y_1},\dots,
\text{Id}_{r \times r} \otimes y_n, \text{Id}_{r \times r}  \otimes \frac{\partial}{\partial y_n}\right)
\end{gather*}
 of
$\text{M}_r(\diff_n)$. Note that $\psi^r_{2n}(c^r_{2n})$ makes sense since $\psi^r_{2n}$ is $G$-basic.

 \begin{proposition} \label{compt}
$f_{\cale,\psi_{2n}}(c_{\cale}(U)) = \psi^r_{2n}(c^r_{2n})$. \end{proposition}

\begin{proof} Note that any element $\mathcal D$ of $\dif(\cale)(U)$
gives a holomorphic $\text{M}_r(\diff_n)$-valued function~$\tilde{D}$ on $J_{\infty}\cale|_{U}$, satisfying $d\tilde{\mathcal
D}+[\Omega,\tilde{\mathcal D}]=0$ where $\Omega$ is the Maurer--Cartan form from Section~3.1.2. It follows that there is a map of complexes
\begin{gather*}
\text{C}_{\bullet}(\dif(\cale)(U)) \rar \text{C}_{\bullet}((\Omega^{\bullet,0}(J_{\infty}\cale|_{U},M_r(\diff_n)),d+[\Omega,-]))\\
\hphantom{\text{C}_{\bullet}(\dif(\cale)(U))}{}  \rar
\text{C}^{\Pi}_{\bullet}(\Omega^{\bullet,0}(J_{\infty}\cale|_{U},M_r(\diff_n))).
\end{gather*}
 The last arrow is by Proposition~\ref{dgatwist}.
Evaluation at $\psi^r_{2n}$ followed by applying the involution multiplying $p$-forms by ${(-1)}^{\lfloor \frac{p}{2} \rfloor}$ therefore
yields a map of complexes
\begin{gather*}
\text{C}_{\bullet}(\dif(\cale)(U)) \rar \Omega^{2n-\bullet,0}(J_{\infty}\cale|_{U}).
\end{gather*}
 Given a section
$\varphi:J_1\cale \rar J_\infty\cale$ and a trivialization of $J_1\cale$ over $U$, one has a composite map
\begin{gather} \label{rewrite}
\text{C}_{\bullet}(\dif(\cale)(U)) \rar \Omega^{2n-\bullet,0}(J_{\infty}\cale|_{U}) \rar \Omega^{2n-\bullet}(J_1\cale|_{U}) \rar
\Omega^{2n-\bullet}(U).
\end{gather} The middle arrow above is pullback by the section $\varphi:J_1\cale \rar J_{\infty}\cale$. The
rightmost arrow above is pullback by the section of $J_1\cale$ arising out of the trivialization of $J_1\cale$ that we have chosen over $U$.
Any section $\varphi: J_1\cale|_{U} \rar J_{\infty}\cale|_{U}$ is {\it unique upto homotopy}. It follows from this that the composite map
\eqref{rewrite} is unique upto homotopy for a f\/ixed trivialization of $J_1\cale$ over $U$. This fact and the proof of Proposition~\ref{freetriv} together imply that the image of $c_{\cale}(U)$ under~\eqref{rewrite} is independent of the precise choice of $\varphi$ and
of trivialization of $J_1\cale$ over $U$. To compute it, we could choose $\varphi$ to be the section taking $\frac{\partial}{\partial z_i}$ to
the formal derivative $\frac{\partial}{\partial y_i}$ and $z_i$ to $y_i+a_i$ at $(a_1,\dots,a_n)$. In this situation, the image of
$c_{\cale}(U)$ is indeed seen to be $\psi^r_{2n}(c^r_{2n})$ in homological degree $2n$.

For any section $\varphi$ of $J_{\infty}\cale|_{U}$ over $J_1\cale$ and for any trivialization of $J_1\cale$ over $U$, the map~\eqref{rewrite} takes $(\mathcal D_0,\dots,\mathcal D_{2n})$ to $\psi^r_{2n}(\hat{\mathcal D_0},\dots,\hat{\mathcal D_n})$ where $\hat{\mathcal
D}$ is the section of $\Omega^{\bullet}(U,\text{M}_r(\diff_n))$ corresponding to $\mathcal D$ (satisfying $d\hat{\mathcal
D}+[\omega,\hat{\mathcal D}]=0$ where $\omega$ is as in Section~3.1.3). Let $\varphi$ be a section of $J_{\infty}\cale|_{U}$ over $J_1\cale$
associated with the Fedosov dif\/ferential $D$ on $\Omega^{\bullet}(\dga{\cale})$ as in Section~3.1.5. Proposition~\ref{proposition8} implies that
$f_{\cale,\psi_{2n}}$ takes $(\mathcal D_0,\dots,\mathcal D_{2n})$ to $\psi^r_{2n}(\hat{\mathcal D_0},\dots,\hat{\mathcal D_n})$. Thus,
$f_{\cale,\psi_{2n}}(c_{\cale}(U)) = \psi^r_{2n}(c^r_{2n})$.
\end{proof}

\textbf{3.3. Properties of $\boldsymbol{f_{\cale,\psi_{2n}}}$.}
Let $\dcat{X}$ denote the bounded derived category of sheaves of $\compl$-vector spaces on
$X$. Then, the map $f_{\cale,\psi_{2n}}$ of Proposition~\ref{mapmade} induces a morphism in~$\dcat{X}$. Recall that the natural embedding of
$\compl$ into $\Omega^{\bullet}(U)$ for any open $U \subset X$ is a~quasi-isomorphism. This induces a quasiisomorphism of complexes of sheaves
$\underline{\compl} \rar \Omega^{\bullet}_X$ on $X$ between the constant sheaf $\underline{\compl}$ and $\Omega^{\bullet}_X$.

\begin{proposition} \label{welldefined} The following diagram commutes in $\dcat{X}$
\begin{gather*}
\begin{CD} \widetilde{\text{\rm hoch}(\dif(\cale))} @>f_{\cale}>>
\Omega^{2n-\bullet}_X\\ @AA{\beta^{-1}_{\cale}}A    @AAA \\ \underline{\compl}[2n] @>>{\frac{[\psi_{2n}]}{[\tau_{2n}]}\text{\rm id}}>
\underline{\compl}[2n]\\ \end{CD}
\end{gather*}
 In particular, the map $f_{\cale}$ in $\dcat{X}$ is independent of the choice of Fedosov connection on~$\dga{\cale}$.
\end{proposition}

\begin{proof} As objects in $\dcat{X}$, $\widetilde{\text{hoch}(\dif(\cale))}$ as well as $\Omega^{2n-\bullet}_X$ are isomorphic to the
shifted constant sheaf $\underline{\compl}[2n]$. Since $\underline{\compl}$ is an injective object in the category of sheaves of
$\compl$-vector spaces on $X$ (since it is f\/labby),
\begin{gather*}
\text{Hom}_{\dcat{X}}(\underline{\compl},\underline{\compl}) \simeq \compl.
\end{gather*}
It follows that it is enough to verify this proposition for each open subset $U$ of $X$ with local holomorphic coordinates on which $\cale$ is
trivial. For such a $U$, the generalized trace map maps~$c_{\cale}(U)$ to a $2n$-cycle mapping to the normalized Hochschild $2n$-cycle
\begin{gather*}
r
\sum_{\sigma \in S_{2n}} \text{sgn}(\sigma)\sigma\left(1,z_1,\frac{\partial}{\partial z_1},\dots,z_n,\frac{\partial}{\partial z_n}\right).
\end{gather*}
 It
follows from Convention~1, Section~2.1
 that $\beta_{\cale}(c_{\cale}(U))=r$. Therefore,
$\frac{1}{r}c_{\cale}(U)$ is a $2n$-cycle mapped to  $1$ by $\beta_{\cale}$. It follows that $f_{\cale}(\beta_{\cale}^{-1}(1))=
\frac{1}{r}f_{\cale}(c_{\cale}(U))= \frac{1}{r}\psi^r_{2n}(c^r_{2n})$. The last equality is by Proposition~\ref{compt}. However $\psi^r_{2n}(c^r_{2n})=
\frac{[\psi^r_{2n}]}{[\tau^r_{2n}]}\tau^r_{2n}(c^r_{2n})$. But $\frac{[\psi^r_{2n}]}{[\tau^r_{2n}]}=\frac{[\psi_{2n}]}{[\tau_{2n}]}$ by
equation~\eqref{eq1}. Further, since  $\tau_{2n}(c_{2n})=1$ (see \cite{fefesh}), $\tau^r_{2n}(c^r_{2n})=r$ by~\eqref{eq1}. Therefore,
$\frac{1}{r}\psi^r_{2n}(c^r_{2n})=\frac{[\psi_{2n}]}{[\tau_{2n}]}$. This proves the desired proposition. \end{proof}

Since $f_{\cale,\psi_{2n}}$ is a map of complexes of sheaves, the map
\begin{gather}
\Gamma_c(X,\widetilde{\text{hoch}(\dif(\cale))} \rar \compl,\qquad
 \label{intcoc} \alpha \mapsto \int_X f_{\cale,\psi_{2n}}(\alpha) \end{gather} descends to a (nonzero) linear functional on
$\text{H}_{0}(\Gamma_c(X,\widetilde{\text{hoch}(\dif(\cale))}))$. We shall denote this linear functional by $\int_X f_{\cale,\psi_{2n}}$.

\textbf{3.4. Proof of Theorem \ref{th}.}  Proposition~\ref{welldefined} implies that $[f_{\cale,\psi_{2n}}(\alpha)]=\frac{[\psi_{2n}]}{[\tau_{2n}]}[\alpha]$ in
$\text{H}^{2n}_c(X,\compl)$ for any $0$-cycle $\alpha$ in $\Gamma_c(X,\widetilde{\text{hoch}(\dif(\cale))}$. Theorem~\ref{th} follows immediately
from this observation and Theorem~\ref{oldrecall2}.

\textbf{3.5. Remarks.}  $f_{\cale,\psi_{2n}}$ induces a map of complexes
\begin{gather*}
\text{C}_{\bullet}(\diffdt(\cale)) \rar \Omega^{2n-\bullet}(X).
 \end{gather*}
 This can equivalently be viewed as a family $\Theta_i  \in \Omega^i(\text{C}^{2n-i}(\diffdt(\cale)))$ of cochain valued forms
satisfying the dif\/ferential equations
\begin{gather*}
d_{\rm DR}\Theta_i = \pm \delta \Theta_{i+1},
\end{gather*}
 where $\delta$ is the dif\/ferential on the Hochschild
cochain complex $\text{C}^{2n-i}(\diffdt(\cale))$. In particular, $\Theta_{2n}$ is a $2n$-form with values in $\text{C}^0(\diffdt(\cale))$.
When $X$ is compact, $\int_X \Theta_{2n}$ is precisely the Hochschild $0$-cocycle~\eqref{intcoc}. This viewpoint seeing~\eqref{intcoc} as
coming from ``integrating $\psi^r_{2n}$ over $X$'' is taken by~\cite{Sh1}.

More generally, there is a modif\/ied cyclic chain complex $\widetilde{\text{Cycl}(\dif(\cale))}$ related closely to
$\widetilde{\text{hoch}(\dif(\cale))}$ (see Section~2.3). The construction of this section can be repeated for a continuous $G$-basic {\it
cyclic} $2n+2p$-cocycle $\nu_{2n+2p}$ (see \cite{pfpota,will}). One obtains a map of complexes
\begin{gather*}
\Gamma_c(X,\widetilde{\text{Cycl}(\dif(\cale))}) \rar \Omega^{2n+2p-\bullet}_c(X)
\end{gather*}
 as a result. When $X$ is compact, and $k \geq p$, the
above map on the $2k$-th homology yields a map
\begin{gather*}
\text{H}^{2n-2k}(X,\compl) \oplus \text{H}^{2n-2k+2}(X,\compl) \oplus \cdots \oplus
\text{H}^{2n}(X,\compl) \rar \text{H}^{2n+2p-2k}(X,\compl)
\end{gather*}
 It would be interesting to understand this map further.

\section[Integrating the lifting formula $\Psi_{2n+1}$]{Integrating the lifting formula $\boldsymbol{\Psi_{2n+1}}$}\label{section4}

 We begin by
outlining the construction of $\Psi_{2n+1}$ in Section~4.1. Section~4.2 is devoted to proving Theorem~\ref{th3}.

\textbf{Notation.} For any algebra $A$, $\text{M}_r(A)$ will denote the algebra of $r \times r$ matrices with entries in~$A$. $\fmat(A)$ will
denote the Lie algebra of f\/inite matrices with entries in $A$. $\text{M}_{\infty}(A)$ will denote the algebra of inf\/inite matrices $M$ with
entries in $A$ such that the $(i,j)$-th entry of $M$ vanishes whenever $|i-j|>C(M)$ where $C(M)$ is a constant (depending on~$M$).

{\bf 4.1.~About the lifting formula $\boldsymbol{\Psi_{2n+1}}$.}

\textbf{4.1.1.} The only explicit fact about the Lifting formula $\Psi_{2n+1}$ that we shall use in this paper is its form. Let $\psdif_n$
denote the algebra $\compl[[y_1,\dots,y_n]][y_1^{-1},\dots,y_n^{-1}][[\partial_1^{-1},\dots,\partial_n^{-1}]][\partial_1,\dots,\partial_n]$
($\partial_i:=\frac{\partial}{\partial y_i}$) of formal pseudodif\/ferential operators. The coef\/f\/icient at
$y_1^{-1}\cdots y_n^{-1}\partial_1^{-1}\cdots \partial_n^{-1}$ def\/ines a linear functional on $\psdif_n$. This functional vanishes on commutators of
elements of $\psdif_n$. We therefore call it a trace and denote it by $\text{Tr}_{\psdif_n}:\psdif_n \rar \compl$ (see~\cite{A}). Recall that
the ``usual'' matrix trace $\text{Tr}_{\fmat}$ yields a linear map from $\fmat(A)$ to $A$ for any algebra $A$. It can then be verif\/ied that
$\text{Tr}_{\psdif_n} \circ \text{Tr}_{\fmat}$ yields a trace on the algebra~$\fmat(\psdif_n)$.

For any $\alpha \in A$, let $\alpha \otimes 1_{\infty}$ denote the diagonal matrix $\alpha \in \text{M}_{\infty}(A)$. Let $D_i$, $1 \leq i
\leq 2n$ be the derivations on $\fmat(\psdif_n)$ given by $D_i=\text{ad}(\text{ln}(x_i) \otimes 1_{\infty})$ for $1 \leq i \leq n$ and
$D_i=\text{ad}(\text{ln}(\partial_{i-n}) \otimes 1_{\infty})$ for $n+1 \leq i \leq 2n$. We recall from~\cite{Sh2} that $[D_i,D_j]=
\text{ad}(Q_{ij} \otimes 1_{\infty})$ for some elements $Q_{ij}$ of $\psdif_n$.

Consider the set $\mathcal S$ of all markings of the interval $[1,2n-1]$ such that
\begin{enumerate}\itemsep=0pt
\item[(i)] Only f\/initely many integral points are marked.
\item[(ii)] The distance between any two distinct marked points is at least $2$.
\end{enumerate}

Note that the ``empty'' marking where no point is marked is also an
element of $\mathcal S$.

Let $t \in S$ be a marking of $[1,2n-1]$ marking the integers $i_1,\dots,i_k$. Def\/ine
\begin{gather*}
\mathcal{O}(t)(A_1,\dots,A_{2n+1})=
\text{Alt}_{A,D}(\text{Tr}_{\psdif_n} \circ \text{Tr}_{\fmat}(P_{1,t} \circ \cdots  \circ P_{2n+1,t})),
\end{gather*}
 where
\begin{gather*}
 \text{if} \ \ j \ \ \text{is marked,} \quad
P_{j,t}=A_jQ_{j,j+1} \quad \text{and} \quad P_{j+1,t}= A_{j+1},\\
P_{j,t}=D_j(A_j) \ \ \text{if} \ \ j \ \ \text{and} \ \ j-1 \ \ \text{are not marked and} \ \ j \neq 2n+1,\\
P_{2n+1,t}=A_{2n+1}.
\end{gather*}
 Note that $0$ and $2n$ should be thought of as unmarked by default in the above formula. If $t$ is the
``empty'' marking,
\begin{gather*}
{\mathcal O}(t)(A_1,\dots,A_{2n+1})= \text{Alt}_{A,D}(\text{Tr}_{\psdif_n} \circ
\text{Tr}_{\fmat}(D_1(A_1)\cdots D_{2n}(A_{2n})A_{2n+1})).
\end{gather*}

\begin{theorem}[\cite{Sh2}]\label{th6} The linear functional
\begin{gather*}
(A_1,\dots,A_{2n+1}) \mapsto \sum_{t \in \mathcal S} {\mathcal O}(t)(A_1,\dots,A_{2n+1})
\end{gather*}
 is a
$2n+1$-cocycle in $\text{\rm C}^{2n+1}_{\text{\rm Lie}}(\fmat(\psdif_n);\compl)$.
\end{theorem}

 The natural
inclusion of algebras $\diff_n \subset \psdif_n$ extends to an inclusion $\fmat(\diff_n) \subset \fmat(\psdif_n)$. The formula from Theorem~\ref{th6}
therefore, gives us the cocycle {\samepage
\begin{gather*}
\Psi_{2n+1} \in \text{C}^{2n+1}_{\text{Lie}}(\fmat(\diff_n),\compl).
\end{gather*}
 $\Psi_{2n+1}$ is referred to as a
lifting formula.}

\textbf{4.1.2.}  Note that for any algebra $A$, we have isomorphisms
\begin{gather*}
\text{M}_m(A) \otimes \text{M}_r(\compl) \simeq \text{M}_{rm}(A)
\end{gather*}
 of
algebras. Taking the direct limit of these isomorphisms, we obtain a map of algebras (and therefore, Lie algebras)
\begin{gather*}
i_m: \ \ \fmat(\text{M}_m(A))
\simeq \text{M}_m(A) \otimes \fmat(\compl) \rar \fmat(A).
\end{gather*}
 It follows that $\Psi_{2n+1}$ yields a cocycle in
$\text{C}^{2n+1}_{\rm lie}(\fmat(\text{M}_m(\diff_n)),\compl)$ as well. We denote this cocycle by~$\Psi^m_{2n+1}$.

For $M \in \text{M}_m(\diff_n)$, $M \otimes \text{id}$ shall refer to the element $i_m(M \otimes \text{id})$ of $\text{M}_{\infty}(\diff_n)$.
For any $p>m$, let $\iota_{m,p}:\text{M}_m(\diff_n) \rar \text{M}_{p}(\diff_n)$
denote the natural embedding obtained by ``padding with
zeros''.

\begin{proposition}\label{proposition16}
 If $A_1,\dots,A_{2n},M \in \text{\rm M}_m(\diff_n)$ and $B_1,\dots,B_{2n} \in \text{\rm M}_p(\diff_n)$ then
\begin{gather*}
\Psi^{m+p}_{2n+1}((A_1\oplus B_1)\otimes \text{\rm id},\dots,(A_{2n}
 \oplus B_{2n}) \otimes \text{\rm id},\iota_{m,m+p}(M) \otimes
 \text{E}_{1,1})\\
 \qquad{} = \Psi^m_{2n+1}(A_1 \otimes \text{\rm id},\dots,A_{2n}
 \otimes \text{\rm id},M \otimes \text{\rm E}_{1,1}).
 \end{gather*}
\end{proposition}

\begin{proof} It is enough to show that for any $t \in \mathcal S$, \begin{gather}  {\mathcal O}(t)((A_1\oplus B_1)\otimes
\text{id},\dots,(A_{2n}
 \oplus B_{2n}) \otimes \text{id},\iota_{m,m+p}(M) \otimes
 \text{E}_{1,1})\nonumber\\
 \qquad{}  = {\mathcal O}(t)(A_1 \otimes \text{id},\dots,A_{2n}
 \otimes \text{id},M \otimes \text{E}_{1,1}).\label{a1}
 \end{gather}
The summands on the left hand side of~\eqref{a1} are of the form $\text{Tr}_{\psdif_n} \circ \text{Tr}_{\fmat}(X_1 \circ \cdots \circ X_{r})$
where $X_i=(A_j \oplus B_j) \otimes \text{id}$ or $X_i= \iota_{m,m+p}(M) \otimes
 \text{E}_{1,1}$ or $X_i = \lambda \otimes \text{id}$ for some
 element $\lambda$ of $\psdif_n$. It is easy to see that this
 summand does not change if each $X_i$ of the form $X_i=(A_j \oplus B_j) \otimes \text{id}$
 is replaced by $A_j \otimes \text{id}$. Doing this however transforms
 the sum on the left hand side to that on the right hand side.
\end{proof}

{\bf 4.2.~Constructing Shoikhet's holomorphic noncommutative residue in general.}

\textbf{4.2.1.} Let $\mathfrak{g}$ be a DG-Lie algebra. Then, $\{C_{p,q}:=(\wedge^p\mathfrak{g})^{-q} | p \geq 1 \}$ becomes a double complex
whose horizontal dif\/ferential is the Chevalley--Eilenberg dif\/fe\-ren\-tial and whose vertical dif\/fe\-ren\-tial is that induced by the dif\/ferential
intrinsic to $\mathfrak{g}$. We denote the complex $\text{Tot}^{\oplus}(C_{\bullet,\bullet})$ by~$\text{C}^{\rm lie}_{\bullet}(\mathfrak{g})$.
Similarly we denote $\text{Tot}^{\Pi}(C_{\bullet,\bullet})$ by $\text{C}^{\Pi,{\rm lie}}_{\bullet}(\mathfrak{g})$. There is a natural map of
complexes from $\text{C}^{\rm lie}_{\bullet}(\mathfrak{g})$ to $\text{C}^{\Pi,{\rm lie}}_{\bullet}(\mathfrak{g})$. Let $\omega \in \mathfrak{g}^1$ be a
Maurer--Cartan element. Let $\mathfrak{g}_{\omega}$ denote the twisted Lie algebra whose underlying dif\/ferential is $d+[\omega,-]$.
\begin{proposition} \label{twist2} There is a natural map of complexes
\begin{gather*}
\text{C}^{\Pi,{\rm lie}}_{\bullet}(\mathfrak{g_{\omega}}) \rar
\text{C}^{\Pi,{\rm lie}}_{\bullet}(\mathfrak{g}), \\
(g_0,\dots,g_k) \mapsto \sum_{j \geq 0} \frac{1}{j!}(g_0,\dots,g_k,\omega,\omega,\dots,\omega) \qquad (j  \ \ \omega\text{'s}) .
\end{gather*}
 \end{proposition}

\begin{proof} Since this proposition is completely analogous to Proposition~2.4 of \cite{enfe}, we shall only sketch the proof. Denote the
dif\/ferential of $\mathfrak{g}$ by $d$. Let $d_{\rm CE}$ denote the Chevalley--Eilenberg dif\/ferential. Let $(\omega)^j:= \omega \wedge \cdots \wedge
\omega$ $j$~times.

{\it Step 1.} One f\/irst notes that
\begin{gather*}
d_{\rm CE}(\omega)^j=\frac{j(j-1)}{2}[\omega,\omega]\wedge (\omega)^{j-2} = -j(j-1)d\omega \wedge
(\omega)^{j-2} = -jd\big(\omega^{j-1}\big).
\end{gather*}
 The middle equality above is because $\omega$ is a Maurer--Cartan element. It follows that if
$\phi_j:=\frac{1}{j!}(\omega)^j$, then
\begin{gather} \label{stair} d_{\rm CE}\phi_j=-d\phi_{j-1}.
\end{gather}

{\it Step 2.} Let $g_0,\dots,g_k$ be homogenous elements of $\mathfrak{g}$. Let $d_i$ denote the degree of $g_i$. Let $\calg G:=g_0 \wedge \cdots
\wedge g_k$. One then verif\/ies (by a direct calculation) that
\begin{gather} \label{stair1} d_{\rm CE}(\calg G \wedge \phi_j)= (d_{\rm CE}\calg
G)\wedge \phi_j +{(-1)}^{k+1}\calg G \wedge d_{\rm CE}\phi_j +{(-1)}^{d_0+\cdots +d_k+k+1}(\text{ad}(\omega)\calg G) \wedge \phi_{k-1}.\!\!\!
\end{gather}

The desired proposition follows from equations~\eqref{stair} and~\eqref{stair1} after inserting the relevant def\/initions and summing over~$k$.
\end{proof}

\textbf{Notation.} In some situations in this section, we f\/ind it better to specify the dif\/ferential of a DG-Lie algebra: if $d$ is the
dif\/ferential on a DG-Lie algebra $\mathfrak{g}$ we often denote $\text{C}^{\rm lie}_{\bullet}(\mathfrak{g})$ and
$\text{C}^{\Pi,{\rm lie}}_{\bullet}(\mathfrak{g})$ by $\text{C}^{\rm lie}_{\bullet}(\mathfrak{g},d)$ and $\text{C}^{\Pi,{\rm lie}}_{\bullet}(\mathfrak{g},d)$
respectively.

\textbf{4.2.2.} Let $D$ be the Fedosov dif\/ferential on $\Omega^{\bullet}(\dga{\cale})$ chosen as in Section~3.1.4. This choice ensures that
there is a morphism $\difdt(\cale) \rar \Omega^{\bullet}(\dga{\cale})$ of sheaves of DG-algebras on $X$. Let $U$ be an open subset of $X$ on
which $J_1\cale$ is trivial. Recall that any trivialization of $J_1\cale$ on an open subset $U$ of $X$ induces an term by term isomorphism of
sheaves DG-algebras $(\Omega^{\bullet}(\dga{\cale})_U,D) \rar (\Omega^{\bullet}_U(\text{M}_r(\diff_n)),d+[\omega,-])$. Also, $\omega \otimes
1_N$ is a Maurer--Cartan element in $\fmat{(\Omega^{\bullet}(U,\text{M}_r(\diff_n)))}$ for any suf\/f\/iciently large natural number~$N$. Here $1_N$ denotes the $N \times N$ identity matrix ``padded with $0$'s'' on the right and bottom to obtain an element of~$\fmat{(\compl)}$. We will continue to denote this element by $\omega$ for notational brevity.

One therefore has the following composite map of complexes
\begin{gather} \label{mainmap} \begin{CD}
\text{C}^{\rm lie}_{\bullet}(\fmat{(\Omega^{\bullet}(U,\dga{\cale}),D)}) @>>>
 \text{C}^{\rm lie}_{\bullet}(\fmat{(\Omega^{\bullet}(U,\text{M}_r(\diff_n)),d+[\omega,-])})\\
@V{\theta}VV @VVV\\ \text{C}^{\Pi,{\rm lie}}_{\bullet}(\fmat{(\Omega^{\bullet}(U,\text{M}_r(\diff_n)),d)}) @<<<
\text{C}^{\Pi,{\rm lie}}_{\bullet}(\fmat{(\Omega^{\bullet}(U,\text{M}_r(\diff_n)),d+[\omega,-])})\\ \end{CD}
 \end{gather} The horizontal arrow on
top is from the isomorphism of $(\Omega^{\bullet}(U,\dga{\cale}),\!D)$ with $(\Omega^{\bullet}(U,\!\text{M}_r(\diff_n)),d+[\omega,-])$ induced by
$\varphi$. The vertical arrow on the right is the natural map mentioned in Section~4.2.1. The horizontal arrow on the bottom is from
Proposition~\ref{twist2}. Let $\Xi_{2n+1}$ be any continuous $2n+1$ cocycle in $\text{C}^{2n+1}_{\rm lie}(\fmat{(\text{M}_r(\diff_n))},\compl)$.
As in Proposition~\ref{constr1}, Section~3.2.3, eva\-lua\-tion at $ \Xi_{2n+1}$ yields a map of complexes from
$\text{C}^{\Pi,{\rm lie}}_{\bullet}(\fmat{(\Omega^{\bullet}(U,\text{M}_r(\diff_n)),d)})[1]$ to $\Omega^{2n-\bullet}(U)$ (with the dif\/ferential on
$\Omega^{2n-p}$ being ${(-1)}^pd_{\rm DR}$). Composing this map with $\theta$, and applying the involution multiplying a $p$-form by
${(-1)}^{\lfloor \frac{p}{2} \rfloor}$, we obtain a map of complexes
\begin{gather*} %\label{localshkt1}
(\varphi,\Xi_{2n+1})_*: \ \
\text{C}^{\rm lie}_{\bullet}(\fmat{(\Omega^{\bullet}(U,\dga{\cale}),D)})[1] \rar \Omega^{2n-\bullet}(U).
\end{gather*}
Explicitly, if
$\mu \in \text{C}^{\rm lie}_{p}(\fmat{(\Omega^{\bullet}(U,\dga{\cale}),D)})$,
\begin{gather} \label{formula1}(\varphi,\Xi_{2n+1})_*(\mu)=
{(-1)}^{\lfloor \frac{2n-p}{2} \rfloor} \sum_k \frac{1}{k!} \Xi_{2n+1}(\hat{\mu} \wedge \omega^k),
\end{gather} where $\hat{\mu}$ is the image
of $\mu$ in $\text{C}^{\rm lie}_{p}(\fmat{(\Omega^{\bullet}(U,\text{M}_r(\diff_n)),d+[\omega,-])})$.  In particular, we obtain a map of
complexes
\begin{gather*}
 \text{C}^{\rm lie}_{\bullet}(\Gamma(U,\fmat(\difdt(\cale))))[1] \rar \Omega^{2n-\bullet}(U) .
 \end{gather*}

\begin{proposition} The above map
extends to a map
\begin{gather*}
\lambda(\Xi_{2n+1}): \ \ \widetilde{\text{\rm Lie}}(\dif(\cale))[1] \rar \Omega^{2n-\bullet}_U
\end{gather*}
 of complexes of sheaves of
$\compl$-vector spaces on $U$. \end{proposition}

\begin{proof} This follows from the continuity of $\Xi_{2n+1}$. The argument proving this is
completely analogous to the proof of Proposition~\ref{mapmade}. However, since $\Xi_{2n+1}$ may not be $G:=GL(n,\compl) \times GL(r,\compl)$-basic, we
can only guarantee the existence of a map of complexes of sheaves over $U$. \end{proof}

\begin{corollary}
If $X$ is complex parallelizable and $\cale$ is trivial, $\lambda(\Xi_{2n+1})$ yields a map of complexes of sheaves
$\widetilde{\text{\rm Lie}}(\dif(\cale))[1] \rar \Omega^{2n-\bullet}_X$ on $X$. \end{corollary}

\begin{proof} The obstruction to globalizing
$\lambda(\Xi_{2n+1})$ comes from the fact that there is no consistent way of choosing a section of $J_1\cale$ over $X$ in general.
\end{proof}

\begin{remark}
%\textbf{Remark.}{\it
If $\Xi_{2n+1}$ is $GL(n) \times GL(r)$-basic, $\lambda(\Xi_{2n+1})$ yields a map of complexes of sheaves $\widetilde{\text{Lie}}(\dif(\cale))[1] \rar \Omega^{2n-\bullet}_X$ for arbitrary (not necessarily complex parallelizable) $X$ and arbitrary (not necessarily trivial) $\cale$. The proof of this assertion is completely analogous to that of Proposition~\ref{freetriv}.
\end{remark}

\textbf{4.2.3. Properties of $\boldsymbol{\lambda(\Xi_{2n+1})}$.}
\begin{proposition} \label{homl}
For a fixed trivialization of $J_1\cale$, the map induced by $\lambda(\Xi_{2n+1})$ in $\dcat{U}$ depends only on the
class of $\Xi_{2n+1}$ in $\text{\rm H}^{2n+1}_{\rm lie}(\text{\rm M}_r(\diff_n))$.
\end{proposition}

\begin{proof} \sloppy We will show that if $\Xi_{2n+1}$ is a coboundary, then $(\varphi,\Xi_{2n+1})_*$ is null-homotopic. The homotopy we shall
provide will automatically yield a homotopy between $\lambda(\Xi_{2n+1})$ and~$0$. Recall that the map $(\varphi,\Xi_{2n+1})_*$ was obtained
by composing the map $\theta$ from the commutative diagram~\eqref{mainmap} with evaluation at $\Xi_{2n+1}$ followed by multiplying $p$-forms
by ${(-1)}^{\lfloor \frac{p}{2} \rfloor}$. If $\Xi_{2n+1}=d\beta$, evaluation at $\beta$ gives a homotopy between evalua\-tion at $\Xi_{2n+1}$
and $0$. This proves that $(\varphi,\Xi_{2n+1})_*$ is null-homotopic when $\Xi_{2n+1}=d\beta$. Explicitly, if $\mu \in
\text{C}^{\rm lie}_{p}(\fmat{(\Omega^{\bullet}(U,\dga{\cale}),D)})$, our homotopy maps $\mu$ to ${(-1)}^{\lfloor \frac{2n-1-p}{2} \rfloor} \sum_k
\frac{1}{k!} \beta(\hat{\mu} \wedge \omega^k)$ where $\hat{\mu}$ is the image of $\mu$ in
$\text{C}^{\rm lie}_{p}(\fmat{(\Omega^{\bullet}(U,\text{M}_r(\diff_n)),d+[\omega,-])})$. It is easy to see that since $\beta$ is a continuous
cochain, our homotopy yields a homotopy between $\lambda(d\beta)$ and the $0$ map.
\end{proof}

\begin{proposition} \label{wd1}
The map induced by $\lambda(\Xi_{2n+1})$ in $\dcat{U}$ is independent of the trivialization of $J_1\cale$ chosen.
\end{proposition}

\begin{proof} \sloppy Let $\mathfrak{h}$ be a Lie subalgebra of $\mathfrak{g}$. Let $\text{H}^{\bullet}_{\rm lie}(\mathfrak{g},\mathfrak{h})$ denote the
Lie algebra cohomology of~$\mathfrak{g}$ relative to $\mathfrak{h}$ (with coef\/f\/icients in the trivial module). The proof of Lemma~5.2 of~\cite{fefesh} goes through with minor modif\/ications to show that
\begin{gather*}
\text{H}^{2n+1}_{\rm lie}(\fmat(\text{M}_r(\diff_n)),\fmat \oplus \mathfrak{gl}_n \oplus
\mathfrak{gl}_r) \simeq \text{H}^{2n+1}_{\rm lie}(\fmat(\text{M}_r(\diff_n))).
\end{gather*}
 It follows that there exists a $G$-basic cocycle
$\Theta_{2n+1}$ such that $[\Theta_{2n+1}]=[\Xi_{2n+1}]$ in $\text{H}^{2n+1}_{\rm lie}(\fmat(\text{M}_r(\diff_n)))$. By Proposition~\ref{homl},
$\lambda(\Xi_{2n+1})=\lambda(\Theta_{2n+1})$ as a morphism in $\dcat{U}$. The proof that $\lambda(\Theta_{2n+1})$ is independent of the chosen
trivialization of $J_1\cale$ (as a morphism in $\dcat{U}$) is identical that of Proposition~\ref{freetriv} (with the obvious modif\/ications).
\end{proof}

\begin{proposition} \label{wd2} The map induced by $\lambda(\Xi_{2n+1})$ in $\dcat{U}$ is independent of the Fedosov differential $D$ on
$\Omega^{\bullet}(\dga{\cale})$ $($for a fixed trivialization of $J_1\cale$ over~$U)$.
\end{proposition}

\begin{proof} From the discussion in Section 2.4.1, $\widetilde{\text{Lie}}(\dif(\cale))[1]$ is isomorphic to $\underline{\compl}[2n] \oplus
V$ in $\dcat{U}$ where $V$ is a direct sum of shifted constant sheaves concentrated in homological degree $\geq 2n+2$. Also,
$\Omega^{2n-\bullet}$ is isomorphic to $\underline{\compl}[2n]$ in $\dcat{U}$. Since $\underline{\compl}$ is an injective object in the
category of sheaves of $\compl$-vector spaces, $\text{Hom}_{\dcat{U}}(V,\underline{\compl}[2n])=0$. It therefore, suf\/f\/ices to show that
$\lambda(\Xi_{2n+1})$ applied to a f\/ixed $2n$ cycle generating the $2n$-st homology of $\Gamma(U,\widetilde{\text{Lie}}(\dif(\cale))[1])$ is
independent of the choice of Fedosov dif\/ferential.

In fact, since $J_1\cale$ is trivial over $U$, the $2n$-th homology of $\Gamma(U,\widetilde{\text{Lie}}(\dif(\cale))[1])$ is in fact generated
by a $2n$-cycle $c_{\cale}^{\rm lie}(U)$ in $\widehat{\text{C}}^{\rm lie}_{\bullet}(\Gamma(U,\fmat(\dif(\cale))))[1]$ (see Section~2.4.1).

As in the proof of Proposition~\ref{compt}, given any ($G$-equivariant) section $\varphi:J_1\cale \rar J_\infty\cale$, one has a composite map of
complexes \begin{gather} \label{keystep} \widehat{\text{C}}^{\rm lie}_{\bullet}(\Gamma(U,\fmat(\dif(\cale))))[1] \rar
\Omega^{2n-\bullet,0}(J_{\infty}\cale|_{U}) \rar \Omega^{2n-\bullet}(J_1\cale) \rar \Omega^{2n-\bullet}(U).
\end{gather} The
``restriction'' of the f\/irst map above to
 $\text{C}^{\rm lie}_{\bullet}(\Gamma(U,\fmat(\dif(\cale))))[1]$ is constructed just like the analogous map in the proof of Proposition~\ref{compt}. The
 middle arrow is $\varphi^*$ where
$\varphi:J_1\cale \rar J_{\infty}\cale$. The third is the pullback by the chosen section of $J_1\cale$. Since $\varphi$ is unique upto
homotopy, the image of $c^{\rm lie}_{\cale}(U)$ under the composite map~\eqref{keystep} is independent of $\varphi$. Finally, Proposition~\ref{proposition8} tells
us that if $\varphi$ is the section associated with the Fedosov connection $D$ on $\dga{\cale}$ (as in Section~3.1.5), the image of
$c^{\rm lie}_{\cale}(U)$ under the composite map \eqref{keystep} coincides with the image of $c^{\rm lie}_{\cale}(U)$ under the ``restriction'' of
$\lambda(\Xi_{2n+1})$ to $\widehat{\text{C}}^{\rm lie}_{\bullet}(\Gamma(U,\fmat(\dif(\cale))))[1]$. This proves the desired proposition.
\end{proof}

Propositions~\ref{wd1} and~\ref{wd2} show that the map induced by $\lambda(\Xi_{2n+1})$ in $\dcat{U}$ is independent of the choices made to
def\/ine it. Let $\calg F$ be another holomorphic vector bundle on $X$. Let $\iota$ denote the natural map of sheaves of DG-algebras between
$\difdt(\cale)$ and $\difdt(\cale \oplus \calg F)$. Note that $\iota$ induces a map of sheaves $\bar{\iota}:
\widetilde{\text{Lie}}(\dif(\cale))[1] \rar \widetilde{\text{Lie}}(\dif(\cale \oplus \calg F))[1]$ on $X$. Let $r$ and $s$ be the ranks of
$\cale$ and $\calg F$ respectively. Recall that for every positive integer $m$ we have the lifting cocycle $\Psi^m_{2n+1}$ described in
Section~4.1.2 ($\Psi^1_{2n+1}=\Psi_{2n+1}$).

\begin{proposition} \label{book2} The following diagram commutes in $\dcat{U}$
\begin{gather*}
\begin{CD} \widetilde{\text{\rm Lie}}(\dif(\cale))[1] @>\bar{\iota}>>
\widetilde{\text{\rm Lie}}(\dif(\cale \oplus \calg F))[1]\\ @VV{\lambda(\Psi^r_{2n})}V   @V{\lambda(\psi^{r+s}_{2n})}VV\\ \Omega^{2n-\bullet}_U
@>{\text{\rm id}}>> \Omega^{2n-\bullet}_U\\ \end{CD}
\end{gather*}
 \end{proposition}

 \begin{proof} Recall that there is a natural map of complexes from
$\widehat{\text{C}}^{\rm lie}_{\bullet}(\Gamma(U,\fmat(\dif(\cale))))[1]$ to $\Gamma(U,\widetilde{\text{Lie}}(\dif(\cale))[1])$ (and similarly for
$\cale \oplus \calg F$). As observed while proving Proposition~\ref{proposition5}, this map is an isomorphism on homology, and the constant sheaf of $U$
corresponding to the homology of $\widehat{\text{C}}^{\rm lie}_{\bullet}(\Gamma(U,\dif(\cale)))[1]$ is isomorphic to
$\widetilde{\text{Lie}}(\dif(\cale))[1]$ in $\dcat{U}$. It therefore, suf\/f\/ices to show that the following diagram commutes upto homology
\begin{gather} \label{a5} \begin{CD} \widehat{\text{C}}^{\rm lie}_{\bullet}(\Gamma(U,\fmat(\dif(\cale))))[1] @>\bar{\iota}>>
\widehat{\text{C}}^{\rm lie}_{\bullet}(\Gamma(U,\fmat(\dif(\cale \oplus \calg F))))[1]\\ @VV{\lambda(\Psi^r_{2n})}V
@V{\lambda(\psi^{r+s}_{2n})}VV\\ \Omega^{2n-\bullet}_U @>{\text{\rm id}}>> \Omega^{2n-\bullet}_U\\ \end{CD} \end{gather}

By arguments paralleling the proofs of Propositions \ref{wd1} and \ref{wd2}, we are free to choose the sections $U \rar J_{\infty}\cale|_{U}$
and $U \rar J_{\infty}(\cale \oplus \calg F)|_{U}$ that we shall use. Fix holomorphic coordinates $z_1,\dots,z_n$ on $U$. Fix a trivialization
of $\cale$ over $U$. Also f\/ix a trivialization of $\calg F$ over $U$. If~$\{e_i\}$ and~$\{f_j\}$ are the ordered bases of $\Gamma(U,\cale)$
and $\Gamma(U,\calg F)$ (over $\Gamma(U,\mathcal O)$) we have chosen in our trivializations, we trivialize $\cale \oplus \calg F$ by choosing
the ordered basis $\{e_1,\dots,e_r,f_1,\dots,f_s\}$ of \mbox{$\Gamma(U,\cale \oplus \calg F)$}. With these choices, we obtain sections of $J_1\cale$ and
$J_1(\cale \oplus \calg F)$ over $U$. The chosen coordinate system and trivialization of $\cale$ identify $\Gamma(U,\dif(\cale))$ with
$\text{M}_r(\diff(U))$. The ($G$-equivariant) section of $J_\infty\cale$ we choose is the one that maps $\frac{\partial}{\partial z_i}$ to
$\frac{\partial}{\partial y_i}$ and $f(z_1,\dots,z_n)$ to
\begin{gather*}
(a_1,\dots,a_n) \mapsto \sum_I \frac{\partial f}{\partial z^I}\Big|_{(a_1,\dots,a_n)}y^I.
\end{gather*}
Let $D \in \Gamma(U,\fmat(\dif(\cale)))$ be arbitrary. Let $\hat{D}$ denote the f\/lat section of
$\text{C}^{\infty}(U,\fmat(\text{M}_r(\diff_n)))$ corresponding to $D$. Note that by our choices of sections of $J_{\infty}\cale$ and
$J_{\infty}(\cale \oplus \calg F)$,
\begin{gather*}
\widehat{\iota(D)}= \iota_{r,r+s}\hat{D}.
\end{gather*} It follows from equation \eqref{formula1} and
Proposition~\ref{proposition16} that the following diagram commutes literally (not just upto homology)
\begin{gather} \label{a6} \begin{CD} {\text{C}}^{\rm lie}_{\bullet}(\Gamma(U,\fmat(\dif(\cale))))[1] @>\bar{\iota}>>
{\text{C}}^{\rm lie}_{\bullet}(\Gamma(U,\fmat(\dif(\cale \oplus \calg F))))[1]\\ @VV{\lambda(\Psi^r_{2n})}V   @V{\lambda(\psi^{r+s}_{2n})}VV\\
\Omega^{2n-\bullet}_U @>{\text{id}}>> \Omega^{2n-\bullet}_U\\ \end{CD}
\end{gather}
Since ${\text{C}}^{\rm lie}_{\bullet}(\Gamma(U,\fmat(\dif(\cale))))[1]$ is dense in
$\widehat{\text{C}}^{\rm lie}_{\bullet}(\Gamma(U,\fmat(\dif(\cale))))[1]$ (and similarly for $\cale \oplus \calg F$), and since all the maps
involved in diagram~\eqref{a6} are continuous, commutativity of diagram~\eqref{a5} follows (note that this ``commutativity on the nose'' of
diagram~\eqref{a5} is with our convenient choice of section of $J_{\infty}\cale|_{U}$ over $J_1\cale|_{U}$ as well as our choice
trivialization of $J_1\cale$ over $U$). This proves the desired proposition.
\end{proof}

We remark that our proof of Propositions~\ref{wd1}, \ref{wd2} and~\ref{book2} go through for arbitrary complex parallelizable
manifolds. Further, if $\Xi_{2n+1}$ is replaced by a $GL(n) \times GL(r)$-basic cocycle having the same class in cohomology, Proposition~\ref{wd1} is replaced by the obvious analog of Proposition~\ref{freetriv}, which goes through for arbitrary $X$ and arbitrary $\cale$ with bounded geometry. Similarly, if $\Psi_{2n+1}$ is replaced by a $GL(n)$-basic cocycle representing the same cohomology class, Propositions \ref{wd2} and \ref{book2} hold for arbitrary $X$ and arbitrary $\cale$ and $\calg F$ with bounded geometry.

\textbf{4.2.4. Proof of Theorem~\ref{th3}.} {\it For this subsection, we shall assume that $X$ is complex parallelizable, $\cale$ is trivial and $\cale$ has bounded
geometry. The first two assumptions can be removed if $\Psi_{2n+1}$ is replaced by a $GL(n)$-basic cocycle representing the same cohomology class.} By Propositions~\ref{wd1} and~\ref{wd2}, for any choice used in the construction of $\lambda(\Psi^r_{2n+1})$, the map
\begin{gather*}
\int_X
\lambda(\Psi^r_{2n+1}): \ \ \Gamma_c(X,\widetilde{\text{Lie}}(\dif(\cale))[1]) \rar \compl, \qquad \alpha \mapsto \int_X
\lambda(\Psi^r_{2n+1})(\alpha)
\end{gather*}
 induces the same map on the $0$-th homology of $\Gamma_c(X,\widetilde{\text{Lie}}(\dif(\cale))[1])$. On the
other hand, we saw in Section~2.4.2 that $\fls^{{\rm lie},\cale}$ induces a map on the $0$-th homology of
$\Gamma_c(X,\widetilde{\text{Lie}}(\dif(\cale))[1])$  as well. By Corollary~\ref{cor1},
$\text{H}_0(\Gamma_c(X,\widetilde{\text{Lie}}(\dif(\cale))[1]))$ is a $1$-dimensional $\compl$-vector space. It follows that as linear
functionals on $\text{H}_0(\Gamma_c(X,\widetilde{\text{Lie}}(\dif(\cale))[1]))$,
\begin{gather*}
\int_X \lambda(\Psi^r_{2n+1})=C(X,\cale).\fls^{{\rm lie},\cale}
\end{gather*}
 Propositions \ref{book1} and \ref{book2} together with the nontriviality of $\fls^{{\rm lie},\cale}$ on homology imply that
$C(X,\cale)=C(X,\cale \oplus \calg F)$ for any vector bundle $\calg F$ with bounded geometry on $X$. This shows that $C(X,\cale)$ is
independent of $\cale$.

{\it For the rest of this subsection} assume that $\cale=\strcc{X}$. Let $U$ be an open subset of $X$ with holomorphic coordinates that
identify $U$ with an open disk in $\compl^n$. Choose a nontrivial $0$-cycle~$\alpha$ of $\Gamma_c(U,\widetilde{\text{Lie}}(\dif(\cale))[1])$.
Note that after making the necessary choices, the construction of $\lambda(\Psi_{2n+1})$ is local in nature. It follows that{\samepage
\begin{gather}\label{localnature1}
\int_U  \lambda(\Psi_{2n+1})(\alpha)=\int_X \lambda(\Psi_{2n+1})(j_*\alpha),
\end{gather} where $j:U \rar
X$ is the natural inclusion.}

 Further, if $\phi$ is an element of $\diff^0(X)$ that is compactly supported on $U$,
\begin{gather} \label{localnature2}
\text{str}\big(\phi e^{-t\Delta_{X}}\big)= \text{str}\big(\phi e^{-t\Delta_{{U}}}\big)
\end{gather} for any $t>0$. In
the right hand side of the above equation, we think of $\phi$ and $e^{-t\Delta_{{U}}}$ as endomorphisms of the space of square integrable
sections of $\dol{\strcc{U}}{}\!\!$. To see this, note that if~$p_t(x,y)$ is the kernel of $e^{-t\Delta_{X}}$,
\begin{gather*}
 \text{str}(\phi
e^{-t\Delta_{X}}) =\int_X \text{str}(\phi p_t(x,x))|dx| =\int_U \text{str}(\phi p_t(x,x))|dx| = \text{str}\big(\phi e^{-t\Delta_{{U}}}\big).
\end{gather*}

The second equality above is because $\phi$ is compactly supported on $U$. The third equality above is because the heat kernel on $U$ is
unique once the choice of Hermitian metric on $\dol{\strcc{U}}{}\!\!$ is f\/ixed (see~\cite{Don}). It follows from equation~\eqref{localnature1}
(in the case when $\cale=\strcc{X}$) and equation \eqref{localnature2} that $C(X,\cale)$ is independent of~$X$ as well. This proves Theorem~\ref{th3}.

\textbf{4.2.5.} {\it For this subsection, let $X$ be an arbitrary compact smooth manifold.} Let ${\mathcal F}_{\bullet}$ be a~complex of sheaves on~$X$
such that each ${\mathcal F}_i$ is a module over the sheaf of smooth functions on~$X$. Suppose also that ${\mathcal F}_{\bullet}$ is in
$\dcat{X}$. Let $\mathfrak{U}:=\{U_i\}$ be a f\/inite good cover of $X$. Consider the double complex $\check{C}(\mathfrak{U},{\mathcal
F}_{\bullet})$ where
\begin{gather*}
\check{C}_{p,q}(\mathfrak{U},{\mathcal F}_{\bullet}) = \oplus_{j_1 <\cdots <j_p} \Gamma_c(U_{j_1} \cap\cdots \cap
U_{j_p},{\mathcal F}_q).
\end{gather*}
Here, $\Gamma_c$ denotes ``sections with compact support''. The vertical dif\/ferential in this double
complex is induced by the dif\/ferential on ${\mathcal F}_{\bullet}$. The horizontal dif\/ferential $\delta$ is given by the formula
\begin{gather*}
(\delta\omega)_{j_1,\dots,j_{p-1}} = \sum_{j_0=1}^N (\omega)_{j_0,\dots,j_{p-1}}.
\end{gather*}
 Here, $(\alpha)_{i_1,\dots,i_k}$ denotes the component
of a chain $\alpha \in \check{C}_{k,l}(\mathfrak{U},{\mathcal F}_{\bullet})$ in $\Gamma_c(U_{i_1} \cap\cdots \cap U_{i_k},{\mathcal F}_l)$. We
also follow the convention that if $i_1 <\cdots   <i_k$, $(\alpha)_{i_{\sigma(1)},\dots,i_{\sigma(k)}} = \text{sgn}(\sigma)(\alpha)_{i_1,\dots,i_k}$ for
any permuta\-tion~$\sigma$ of $1,\dots,k$. The following proposition is a trivial modif\/ication of Proposi\-tion~12.12 of Bott and Tu~\cite{bttu}. We therefore, omit its proof.
\begin{proposition} \label{cech} \qquad

$(1)$ The map
\begin{gather*}
\text{\rm Tot}_{\bullet}(\check{C}(\mathfrak{U},{\mathcal F}_{\bullet})) \rar \Gamma_c(X,{\mathcal
F}_{\bullet}), \qquad \alpha \mapsto \sum_i (\alpha)_i
\end{gather*} is a quasi-isomorphism.

$(2)$
$\text{\rm H}_p(\text{\rm Tot}_{\bullet}(\check{C}(\mathfrak{U},{\mathcal F}_{\bullet}))) \simeq {\mathbb H}^{-p}_c(X,{\mathcal F}_{\bullet})$. \end{proposition}

The subscript `$c$' on the right hand side in this proposition stands for ``compact support''. The chain complex ${\mathcal
F}_{\bullet}$ needs to be converted into a cochain complex by inverting degrees in order to make sense of the hypercohomology in the above
proposition.

Let $\alpha$ be $0$-cycle in $\Gamma_c(X,{\mathcal F}_{\bullet})$.

\begin{proposition} \label{split} Suppose that ${\calg F}_{\bullet}$ is acyclic in negative degrees. If $X$ is compact, there exist $0$-cycles
$\alpha_i$ of $\Gamma_c(U_i,{\mathcal F}_{\bullet})$ such that $\sum_i [\alpha_i]=[\alpha]$ in $\text{\rm H}_0(\Gamma_c(X,{\mathcal
F}_{\bullet}))$. \end{proposition}

\begin{proof}\looseness=1 \sloppy We use a ``staircase argument''. By Proposition~\ref{cech}, there exists a $0$-cycle $\hat{\alpha}$ of
$\text{Tot}_{\bullet}(\check{C}(\mathfrak{U},{\mathcal F}_{\bullet}))$ such that $[\hat{\alpha}]$ maps to $[\alpha]$ under the
quasi-isomorphism in Proposition~\ref{cech}. Note that $\hat{\alpha}= \sum_{k \geq 1} \hat{\alpha}_k$ where $\hat{\alpha}_j$ is a
$-k+1$-chain in $\oplus_{j_1 <\cdots <j_k} \Gamma_c(U_{j_1} \cap\cdots \cap U_{j_k},{\mathcal F}_{\bullet})$. Let $N(\hat{\alpha})$ be the largest
integer such that $\hat{\alpha}_{N(\hat{\alpha})} \neq 0$. Then, $\hat{\alpha}_{N(\hat{\alpha})}$ is a cycle in $\oplus_{j_1 <\cdots <j_k}
\Gamma_c(U_{j_1} \cap\cdots \cap U_{j_{N(\hat{\alpha})}},{\mathcal F}_{\bullet})$. Since, ${\mathcal F}_{\bullet}$ is acyclic in negative degrees,
$\hat{\alpha}_{N(\hat{\alpha})}=\delta \beta$ for some $\beta$ where $\delta$ is the dif\/ferential on $\oplus_{j_1 <\cdots <j_k} \Gamma_c(U_{j_1}
\cap\cdots \cap U_{j_{N(\hat{\alpha})}},{\mathcal F}_{\bullet})$. Let $d$ denote the dif\/ferential in the bicomplex
$\text{Tot}_{\bullet}(\check{C}(\mathfrak{U},{\mathcal F}_{\bullet}))$. View $\beta$ as a chain in
$\text{Tot}_{\bullet}(\check{C}(\mathfrak{U},{\mathcal F}_{\bullet}))$. Then, $\hat{\alpha}+d\beta$ is homologous to~$\hat{\alpha}$. On the
other hand, $N(\hat{\alpha}+d\beta)=N(\hat{\alpha})-1$. It follows that repeating this procedure eventually gives a $0$-cycle in $\oplus_i
\Gamma_c(U_i,{\mathcal F}_{\bullet})$ homologous to $\hat{\alpha}$. This proves the desired proposition. \end{proof}

\textbf{4.2.6. Shoikhet's holomorphic noncommutative residue: def\/inition.} {\it Let $X$ be an arbitrary compact complex manifold and let
$\cale$ be an arbitrary vector bundle on $X$.} Let $\mathfrak{U}:=\{U_i\}$ be a f\/inite good cover of $X$. Let $\mathcal{D}$ be a global
holomorphic dif\/ferential operator on $\cale$. Note that $\text{E}_{11}(\mathcal D)$ is a $0$-cycle in
$\Gamma_c(X,\widetilde{\text{Lie}}(\dif(\cale))[1])$. Applying Proposition~\ref{split} to the complex $\widetilde{\text{Lie}}(\dif(\cale))[1]$
of sheaves on $X$, we see that there exist $0$-cycles $\alpha_i$ in $\Gamma_c(U_i,\widetilde{\text{Lie}}(\dif(\cale))[1])$ such that $\sum_i
[\alpha_i]=[\text{E}_{11}(\mathcal D)]$ in $\text{H}_0(\Gamma_c(X,\widetilde{\text{Lie}}(\dif(\cale))[1]))$. Def\/ine the {\it holomorphic
noncommutative residue} of $\mathcal D$ to be the sum
\begin{gather*}
\text{NC}(\mathcal D):= \sum_i \int_{U_i} \lambda(\Psi^r_{2n+1})(\alpha_i).
\end{gather*}
One can make the choices used to def\/ine $\lambda(\Psi^r_{2n+1})$ in the $U_i$ arbitrarily. We also remark that $\lambda(\Psi^r_{2n+1})$ is not
globally def\/ined on $X$ in a direct way (at least as far as we can see). By Theorem~\ref{th3},
\begin{gather*}
 \text{NC}(\mathcal D)=  \sum_i \int_{U_i}
\lambda(\Psi^r_{2n+1})(\alpha_i) = C.\sum_i \fls^{{\rm lie},\cale}(\alpha_i) = C.\fls^{{\rm lie},\cale}\bigg(\sum_i \alpha_i\bigg) \\
\phantom{\text{NC}(\mathcal D)}{} =
C.\fls^{{\rm lie},\cale}(\text{E}_{11}(\mathcal D)) = C.\lim_{t \rar \infty} \text{str}\big({\mathcal D}e^{-t\Delta_{\cale}}\big) =
C.\text{str}(\mathcal D).
\end{gather*}
 This proves Conjecture~3.3 of \cite{Sh1} in greater generality.

\looseness=1
Alternatively, one can sidestep the complications that arose while def\/ining $\text{NC}(\mathcal D)$ as follows. Let $\Theta_{2n+1}$ be a $GL(n)$-basic cocycle representing the cohomology class $[\Psi_{2n+1}]$ (its existence is proven following the proof of Lemma~5.2 in [10]). Construct the cocycles $\Theta^r_{2n+1}$ from $\Theta_{2n+1}$ exactly as the $\Psi^r_{2n+1}$ were constructed from $\Psi_{2n+1}$. The cocycles $\Theta^r_{2n+1}$ are $GL(n) \times GL(r)$-basic representatives of the class $[\Psi^r_{2n+1}]$. Then, $\lambda(\Theta^r_{2n+1})$ is a globally def\/ined $2n$-form and def\/ine
\begin{gather*}
\text{NC}(\mathcal D):=\int_X \lambda(\Theta^r_{2n+1})(\text{E}_{11}(\mathcal D)).
\end{gather*}
The two def\/initions of $\text{NC}(\mathcal D)$ given here are equivalent by Proposition~\ref{homl}.

\subsection*{Acknowledgements} I am grateful to Giovanni Felder and Thomas Willwacher for some very useful discussions. This work would not have
reached its current form without their pointing out important shortcomings in earlier versions. I am also grateful to Boris Shoikhet for
useful discussions. I thank the referees of this article for their constructive suggestions. This work was done (prior to my joining my current position) partly at Cornell University and partly at IHES. I am grateful to both these institutions for providing me with a congenial work atmosphere.

\newpage

\pdfbookmark[1]{References}{ref}
\LastPageEnding

\end{document}